%% file: main.tex
\documentclass[twoside]{article}

\usepackage[accepted]{aistats2025}
\usepackage[round]{natbib}

\usepackage{microtype}
\usepackage{graphicx}
\usepackage{subfigure}
\usepackage{booktabs} % for professional tables
\newcommand{\comment}[1]{}
\usepackage{hyperref}
\usepackage{enumitem}% http://ctan.org/pkg/enumitem
\setlist[itemize]{noitemsep, topsep=0pt}

\usepackage{amsmath}
\usepackage{amssymb}
\usepackage{mathtools}
\usepackage{amsthm}
\usepackage{mdframed}
\usepackage[capitalize,noabbrev]{cleveref}

\usepackage{caption}
\usepackage{multirow, hhline}

\usepackage{color}
\usepackage{pifont}
\definecolor{mydarkgreen}{RGB}{39,130,67}

\definecolor{mydarkred}{RGB}{192,47,25}

\theoremstyle{plain}
\newtheorem{theorem}{Theorem}[section]
\newtheorem{proposition}[theorem]{Proposition}
\newtheorem{lemma}[theorem]{Lemma}
\newtheorem{corollary}[theorem]{Corollary}
\newtheorem{example}[theorem]{Example}
\newmdtheoremenv{framedtheorem}[theorem]{Theorem}
\newmdtheoremenv{framedcorollary}[theorem]{Corollary}
\newmdtheoremenv{framedlemma}[theorem]{Lemma}
\newmdtheoremenv{framedexample}[theorem]{Example}
\newmdtheoremenv{framedassumption}[theorem]{Assumption}
\newmdtheoremenv{framedproposition}[theorem]{Proposition}

\theoremstyle{definition}
\newtheorem{definition}[theorem]{Definition}
\newtheorem{assumption}[theorem]{Assumption}
\theoremstyle{remark}
\newtheorem{remark}[theorem]{Remark}

\usepackage[textsize=tiny]{todonotes}

\input{math_commands.tex}

\usepackage{hyperref}
\usepackage{url}
\usepackage{enumitem}
\usepackage[font=small]{caption}

\usepackage{amsthm}
\usepackage{graphicx}

\usepackage{algorithm}
\usepackage[noend]{algpseudocode}

\usepackage{caption}

\algdef{SE}[DOWHILE]{Do}{DoUntil}{\algorithmicdo}[1]{\algorithmicuntil\ #1}%
\begin{document}

% If your paper is accepted and the title of your paper is very long,
% the style will print as headings an error message. Use the following
% command to supply a shorter title of your paper so that it can be
% used as headings.
%
%\runningtitle{I use this title instead because the last one was very long}

% If your paper is accepted and the number of authors is large, the
% style will print as headings an error message. Use the following
% command to supply a shorter version of the authors names so that
% they can be used as headings (for example, use only the surnames)
%
%\runningauthor{Surname 1, Surname 2, Surname 3, ...., Surname n}

\twocolumn[

\aistatstitle{
Improving Stochastic Cubic Newton with Momentum
%Instructions for Paper Submissions to AISTATS 2025
}

\aistatsauthor{ El Mahdi Chayti \And Nikita Doikov \And  Martin Jaggi }

\aistatsaddress{ Machine Learning and Optimization Laboratory (MLO), EPFL } ]

\begin{abstract}
  We study stochastic second-order methods for solving general non-convex optimization problems. We propose using a special version of momentum to stabilize the stochastic gradient and Hessian estimates in Newton's method. We show that momentum provably improves the variance of stochastic estimates and allows the method to converge for any noise level. Using the cubic regularization technique, we prove a global convergence rate for our method on general non-convex problems to a second-order stationary point, even when using only a single stochastic data sample per iteration. This starkly contrasts with all existing stochastic second-order methods for non-convex problems, which typically require large batches. Therefore, we are the first to demonstrate global convergence for batches of arbitrary size in the non-convex case for the Stochastic Cubic Newton. Additionally, we show improved speed on convex stochastic problems for our regularized Newton methods with momentum.
\end{abstract}

\section{Introduction}

%%%%%%%%%%%%%%%%%%%%%%%%%%%%%%%%%%%%%%%%%%%%%%%%%%
%%%%%%%%%%%%%%%%%%%%%%%%%%%%%%%%%%%%%%%%%%%%%%%%%%
\subsection{Motivation}
%%%%%%%%%%%%%%%%%%%%%%%%%%%%%%%%%%%%%%%%%%%%%%%%%%
%%%%%%%%%%%%%%%%%%%%%%%%%%%%%%%%%%%%%%%%%%%%%%%%%%

We are interested in solving general stochastic optimization problems of the form
\begin{equation}
    \label{MainProblem}
    \min_{\vx\in\R^d} f(\vx) := \E_{\xi\sim \mathcal{P}}\Big[f_\xi(\vx)\Big]\,,
\end{equation}
where the functions $f_\xi$ are twice differentiable, not necessarily convex, and only reasonably have access to the distribution $\mathcal{P}$ through its samples.
Problems of this type are important in machine learning applications and can be challenging to solve due to both the \textit{stochasticity} and \textit{ill-conditioning} of $f$.

In this work, we study the class of second-order methods,
which employ a \textit{second-order stochastic oracle} of~\eqref{MainProblem}:
for every $\vx \in \R^d$ and random sample $\xi \in \mathcal{P}$,
we assume to have access to stochastic gradient $\nabla f_{\xi}(\vx)$
and stochastic Hessian matrix $\nabla^2 f_{\xi}(\vx)$.
Incorporating second-order information
into an optimization scheme is known to be a powerful approach
to deal with ill-conditioning of the problem~\citep{polyak1987introduction,nocedal2006numerical,nesterov2018lectures}.

Indeed, for a \textit{full-batch version} of the cubically regularized Newton's method,
which uses the full gradient $\nabla f(\vx)$ and the full Hessian $\nabla^2 f(\vx)$,
it is possible to establish the state-of-the-art global convergence rates
to a second-order stationary point, for general non-convex functions,
and \textit{the rate is provably better} than that of first order 
gradient methods~\citep{nesterov2006cubic}.
Thus, we can see a strong benefit in utilizing second-order information to address ill-conditioning of the problem.

However, in modern large-scale applications, it is often unfeasible to 
have access to a full exact oracle of $f$.
In the presence of stochasticity, such as in machine learning applications, it becomes especially challenging to use noisy gradients and Hessians in optimization algorithms for solving~\eqref{MainProblem}.
While there exist several approaches that use highly accurate
approximations of the gradient and Hessian with subsampling over 
a sufficiently large batch, it remained an important open problem:
\begin{center}
    \textit{Can we design a second-order method that converges globally, for general non-convex functions, for batches as small as one sample?}
\end{center}

\begin{table*}[h!]
	\centering
	\small
	\renewcommand{\arraystretch}{1.3}
        \setlength\tabcolsep{5pt}
	\begin{tabular}{c|c|c|c}
		\hline
		\textbf{Algorithm} & 
		\begin{tabular}{c}
            \textbf{Global} \\
            \textbf{complexity}
            \end{tabular}
            & 
            \begin{tabular}{c}
            \textbf{Arbitrary} \\
            \textbf{accuracy \boldsymbol{$\varepsilon$}}
            \end{tabular}
            &
		  \begin{tabular}{c}
            \textbf{Smooth.} \\
            \textbf{assump.}
            \end{tabular} \\
		\hline
            \begin{tabular}{c}
            Stochastic Gradient Descent (SGD)\\
            \citep{lan2020first}
            \end{tabular}
            & 
            \begin{tabular}{c}
		$F_0 \cdot \cO\Bigl( \frac{{\color{mydarkred}
            \boldsymbol{L_g}}}{\varepsilon^2} 
                + \frac{{\color{mydarkred}
            \boldsymbol{L_g}^{1/2}} \sigma_g}{F_0^{1/2 }\varepsilon^4}  \Bigr)$ 
            \vspace*{2pt}
            \end{tabular}
            & 
            { \color{mydarkgreen}
            \textbf{yes}
            }
            &
            Lipschitz $\nabla f$
		%\cite{ghadimi2013stochastic}  
            \\
            \hline
            \begin{tabular}{c}
            Normalized SGD with momentum\\
            \citep{cutkosky2020momentum}
            \end{tabular} & 
            \begin{tabular}{c}
            $F_0 \cdot \cO\Bigl( 
            \frac{{\color{mydarkred}
            \boldsymbol{L_g}}}{\varepsilon^2}
            + \frac{\sigma_g^{13/3}}{F_0^{7/3} L^{2/3} \varepsilon^{7/3}}
            + \frac{L^{1/2} \sigma_g^{2}}{\varepsilon^{7/2}}
            \Bigr)$
            \vspace*{2pt}
            \end{tabular}
            &
            { \color{mydarkgreen}
            \textbf{yes}
            }
            & 
            Lip.
            $\nabla f$ and $\nabla^2 f$
            
            %\cite{cutkosky2020momentum}
            \\
            \hline
            \begin{tabular}{c}
            Stochastic Cubic Newton (SCN)\\
            \citep{chayti2024unified}
            \end{tabular}& 
            \begin{tabular}{c}
            $
            F_0 \cdot 
            \cO\Bigl( 
            {\color{mydarkgreen}
            \boldsymbol{ \frac{L^{1/2}}{\varepsilon^{3/2}}
            }}
            + \frac{\sigma_h}{\varepsilon^2}
            \Bigr)$
            \vspace*{2pt}
            \end{tabular}
            &
            \!\!\!
            {\color{mydarkred}
            \begin{tabular}{c}
            \textbf{only for} \\
            \boldsymbol{$
            \;\;\; \varepsilon \geq \sigma_g^{3/2}$}
            \end{tabular}
            }
            \!\!\!
            & 
            Lipschitz $\nabla^2 f$
            %\cite{chayti2024unified} 
            \\
            \hline 
            \begin{tabular}{c}
            SCN with momentum \\
            \textbf{(ours)}
            \end{tabular}
            & 
            \begin{tabular}{c}
            $F_0 \cdot \cO\Bigl( 
            {\color{mydarkgreen}
            \boldsymbol{ 
            \frac{L^{1/2}}{\varepsilon^{3/2}}} }
            + 
            {\color{mydarkgreen}
            \boldsymbol{ 
            \frac{L^{1/4} \sigma_h^{1/2}}{\varepsilon^{7/4}}
            }}
            + 
            {\color{mydarkgreen}
            \boldsymbol{ 
            \frac{L^{1/2} \sigma_g^2}{\varepsilon^{7/2}}
            }}
            \Bigr)$
            \vspace*{2pt}
            \end{tabular}
            & 
             { \color{mydarkgreen}
            \textbf{yes}
            }
            & 
            Lipschitz $\nabla^2 f$
            %\textbf{this work}
            \\
		\hline
		
	\end{tabular}
	\caption{A comparison of 
        first-order and second-order methods
        for non-convex optimization
        that use a single stochastic sample
        of gradient and Hessian per iteration.
        Global complexity means the number of iterations (single data samples) required
        to obtain an $\varepsilon > 0$
        accuracy for the expected gradient norm:
        $\| \nabla f(\bar{\vx}) \| \leq \varepsilon$.
        $F_0 := f(\vx_0) - f^{\star}$
        is the initial functional residual,
        $L_g$ denotes the Lipschitz constant
        of the gradient, and $L$ denotes the Lipschitz constant of the Hessian,
        $\sigma_g$ denotes the variance of
        gradient samples and $\sigma_H$
        stands for the variance of the Hessian samples (see precise definitions in Section~\ref{SectionStochCN}).
        Note that $L$ is not affected by adding to our objective $f$ an arbitrary quadratic function, while  $L_g$, used in first-order methods, can become arbitrarily large.
        }
	\label{TableComplexities}
\end{table*}

%%%%%%%%%%%%%%%%%%%%%%%%%%%%%%%%%%%%%%%%%%%%%%%%%%
%%%%%%%%%%%%%%%%%%%%%%%%%%%%%%%%%%%%%%%%%%%%%%%%%%
\subsection{Contributions}
%%%%%%%%%%%%%%%%%%%%%%%%%%%%%%%%%%%%%%%%%%%%%%%%%%
%%%%%%%%%%%%%%%%%%%%%%%%%%%%%%%%%%%%%%%%%%%%%%%%%%

In this paper, we provide a positive answer to this question by designing
an efficient stochastic second-order optimization algorithm
that can work with stochastic gradients and Hessians of \textit{arbitrarily large variance},
converging to a second-order stationary point for non-convex objectives.
The main idea behind our algorithm is to use a special form of \textit{momentum} for stochastic estimates of the gradient and Hessian. Thus, at iteration $t \geq 1$, we use a recursively defined vector
\beq \label{GradUpdateIntro}
\ba{rcl}
\!\!\!\!\!\!
\vg_t & \!\!\!\!=\!\!\!\! & (1 - \alpha) \vg_{t - 1} + \alpha \nabla f_{\xi_t}\bigl( \vx_t + \frac{1 - \alpha}{\alpha}( \vx_t - \vx_{t - 1} ) \bigr)
\ea
\eeq
and matrix
\beq \label{HessUpdateIntro}
\ba{rcl}
\mH_t & = & (1 - \beta) \mH_{t - 1} + \beta \nabla^2 f_{\xi_t}(\vx_t),
\ea
\eeq
where $( \vx_t )_{t \geq 0}$ is a sequence of points generated by the method,
$0 < \alpha, \beta \leq 1$ are momentum parameters,
and $\xi_t$ is a random sample of data.
We start these updates with $\vg_0 = \nabla f_{\xi_0}(\vx_0)$ and $\mH_0 = \nabla f_{\xi_0}(\vx_0)$ for some initial starting point.
Then, to compute the next point $\vx_{t + 1}$, we perform 
the Cubic Newton step from the previous one using freshly estimated $\vg_t$ and $\mH_t$.
Therefore, in a non-degenerate case, one step of our second-order method
can be represented in the following standard form:
\beq \label{CubicNewtonStep}
\ba{rcl}
\vx_{t + 1} & = & \vx_t - \bigl( \mH_t + \gamma_t \mI \bigr)^{-1} \vg_t,
\ea
\eeq
where $\gamma_t > 0$ is a parameter that comes from the cubic regularization
of the second-order model, and which ensures that $\mH_t + \gamma_t \mI \succeq 0$
(see Algorithm~\ref{Mainalg} in Section~\ref{SectionMainTheory} for a formal description of our method).

Note that if we choose $\alpha = 1$ and $\beta = 1$ in our update of $\vg_t$ and $\mH_t$,
we obtain simply $\vg_t = \nabla f_{\xi_t}(\vx_t)$ and $\mH = \nabla^2 f_{\xi_t}(\vx_t)$
which corresponds to a straightforward 'naive' implementation of the stochastic Cubic Newton.
In our work, we show that in the presence of noise, it is essential to pick $\alpha, \beta < 1$,
which adds momentum to the updates and which ensures convergence of the method
for samples with arbitrarily large variance.

We summarize our contributions as follows.

\begin{itemize}
    \item We propose combining the special form of momentum~\eqref{GradUpdateIntro}, \eqref{HessUpdateIntro} with stochastic second-order methods and demonstrate how it can reduce the variance of stochastic samples by simulating large batches using past estimates of the gradient and Hessian.

    \item We prove that the Stochastic Cubic Newton method with momentum enjoys an improved convergence rate, for general non-convex problems (see Table~\ref{TableComplexities}). In particular, we show that the method converges for any choice of batch size, 
    \textit{even just a single sample}. Importantly, our result is obtained under the standard assumptions commonly used in previous literature. 

    \item We also show an improved convergence rate for our method in the convex case (see Section~\ref{SectionConvex}), again for any sample size.
    Notably, we prove this rate for the primal Cubic Newton scheme, without the need
    to use additional extragradient steps.

    \item Illustrative numerical experiments are provided.
\end{itemize}

%\begin{itemize}
%    \item We propose combining momentum with second-order methods that need large batches of data and show how it can be used to simulate large batches from past estimates (of gradient and Hessian).
%    \item We show that for a specific choice of momentum, the Cubic Newton method enjoys an improved convergence rate in the non-convex case and in particular it converges for any choice of batch size, even just a single sample. Most importantly, our result is obtained under common assumptions in previous literature. 
%    \item We also show an improved convergence rate in the convex case without using acceleration, again for any sample size.
%\end{itemize}

%%%%%%%%%%%%%%%%%%%%%%%%%%%%%%%%%%%%%%%%%%%%%%%%%%
%%%%%%%%%%%%%%%%%%%%%%%%%%%%%%%%%%%%%%%%%%%%%%%%%%
\subsection{Related work}
%%%%%%%%%%%%%%%%%%%%%%%%%%%%%%%%%%%%%%%%%%%%%%%%%%
%%%%%%%%%%%%%%%%%%%%%%%%%%%%%%%%%%%%%%%%%%%%%%%%%%

When the objective $f$ is non-convex, problem~\eqref{MainProblem} is especially complicated since finding a global minimum is NP-hard in general \citep{NpH}. 
Hence, the reasonable goal is to look for approximate solutions. \textit{First-order gradient methods} \citep{cauchy1847methode,SGD51,SGD52,lan2020first} only use first-order information (gradients) and can, generally, only guarantee convergence to approximate stationary points i.e., a point with a small gradient norm. For non-convex functions, a stationary point may be a saddle point or even a local maximum, which might not be desirable. To obtain faster convergence rates for ill-conditioned problems and
ensure a second-order stationary condition for an approximate stationary point, we can use second-order information. %(Hessian matrix).

\textbf{Second-order optimization.}
Among many variants of the Newton algorithm,
one of the most theoretically established
is the Cubic Newton method~\citep{griewank1981modification,nesterov2006cubic,cartis2011adaptive1,cartis2011adaptive2},
which enjoys many favorable properties such as faster global convergence (compared to first-order methods), and local quadratic convergence, as for the classical Newton's method~\citep{fine1916newton,bennett1916newton,kantorovich1948newton,polyak2007newton}.

Inexact versions of the Cubic Newton
with approximated gradients and Hessians 
studied intensively in the literature~\citep{cartis2011adaptive1,cartis2011adaptive2,wang2018note},
including 
\textit{low-rank quasi-Newton approximations}~\citep{benson2018cubic,kamzolov2023cubic,scieur2024adaptive,doikov2024spectral},
\textit{stochastic subspace methods}~\citep{hanzely2020stochastic,zhao2024cubic},
and \textit{finite-difference approximations}~\citep{cartis2012oracle,grapiglia2022cubic,doikov2023first}.
However, most of these works assume access to deterministic (full batch) 
oracle information about the objective function, which can be very expensive
or even unfeasible in the case of large datasets.

The main approach to overcome this issue consists in using \textit{stochastic} gradient and Hessian estimates \textit{subsampled} from the data \citep{Subsampled3,Subsampled1,Subsampled2,SCNTrip,ghadimi2017second,cartis2018global,doikov2020convex,TensorCN}. In the subsampling technique, we rely on the use of \textit{sufficiently large batches} of data to reduce the effects of noise on the method's convergence, thereby making the gradient and Hessian estimates highly accurate.
However, as was also noticed by \citep{chayti2024unified}, this technique \textit{does not guarantee convergence} of the Cubic Newton for small batches. This is
in contrast to stochastic first-order methods,
that often possess convergence, starting from an arbitrary initialization and even for one sample per iteration~\citep{lan2020first}. 

In fact, all existing analyses of the subsampled stochastic Cubic Newton method can guarantee convergence only to a ball of radius $\cO(\sigma_g^{3/2})$, where $\sigma_g$ is the variance of the gradient noise defined in the usual manner (see Table~\ref{TableComplexities}). Therefore, if the noise level is large, the method can't guarantee convergence to an arbitrary level of accuracy $\varepsilon$.

\textbf{Convex optimization.} This problem seems to be widely present even in the convex case. For example, \citep{TensorCN} shows convergence in terms of the global functional residual up to a ball $\cO(\sigma_g D)$ for stochastic high-order methods.
At the same time, for convex problems, it is possible
to incorporate stochastic Newton steps
into \textit{extragradient}~\citep{korpelevich1976extragradient,nemirovski2004prox} and
\textit{accelerated}~\citep{Nest83,nesterov2008accelerating}
schemes, which provides
global convergence guarantees for arbitrary sample size~\citep{antonakopoulos2022extra,agafonov2024advancing}.
However, extending this approach to non-convex functions is challenging even for first-order methods. 

In our work, we are mainly focused on non-convex functions; thus we propose a different strategy that consists of reusing past gradient and Hessian estimates to simulate a large batch without ever having to actually sample one. 

\textbf{Momentum.}
The idea of using momentum has been extensively explored for first-order optimization.
It was initially proposed by~\citep{POLYAK19641}
in his Heavy-Ball method,
aimed at improving the convergence properties of classical gradient descent. It was shown in~\citep{POLYAK19641} (see also~\citep{polyak1987introduction})
that the Heavy-Ball method achieves 
the accelerated optimal convergence rate
for deterministic minimization of convex quadratic functions.
Recent results show the benefits of using momentum for deterministic optimization, taking into account the intrinsic distribution of the Hessian spectrum~\citep{scieur2020universal},
H\"older continuity of the Hessian~\citep{marumo2024universal},
and also demonstrate the impossibility of acceleration for general strongly convex functions with Lipschitz gradient~\citep{goujaud2023provable}.

In stochastic optimization,
the use of momentum
in SGD and its variants such as Adam~\citep{kingma2017adammethodstochasticoptimization}
is very popular. Different forms of momentum were established in~\citep{tran2022better,cutkosky2020momentumbasedvariancereductionnonconvex,arnold2019reducingvarianceonlineoptimization} that lead to improved convergence rates. Surprisingly, this simple, yet effective idea, has not been explored for stochastic second-order methods, in particular, the Cubic Newton. The only works we were able to find combining both momentum and the Cubic Newton are \citep{wang2019cubicregularizationmomentumnonconvex,gao2022momentumacceleratedadaptivecubic}, and they both consider a different version of momentum than ours. Most importantly, they consider deterministic optimization without stochasticity.

%\textbf{Challenges of analyzing momentum.} 
Reusing past samples is a very well-used idea in machine learning to deal with stochasticity and smooth out estimates (mostly of gradients as this was mainly used for first-order methods). It is very easy to see that momentum, in this sense, introduces bias to the estimates (since the present is not necessarily the same as the past), quantifying this bias is necessary for analyzing any type of momentum. The variants of momentum in \citep{tran2022better,cutkosky2020momentumbasedvariancereductionnonconvex,arnold2019reducingvarianceonlineoptimization}, are but different ways to algorithmically reduce this bias (under different assumptions). Extending the analyses of momentum to Cubic Newton is particularly challenging as there is no closed-form solution of the update at each step, rendering controlling the bias incurred from using momentum very complicated.

%%%%%%%%%%%%%%%%%%%%%%%%%%%%%%%%%%%%%%%%%%%%%%%%%%
%%%%%%%%%%%%%%%%%%%%%%%%%%%%%%%%%%%%%%%%%%%%%%%%%%
\section{Preliminaries}
%%%%%%%%%%%%%%%%%%%%%%%%%%%%%%%%%%%%%%%%%%%%%%%%%%
%%%%%%%%%%%%%%%%%%%%%%%%%%%%%%%%%%%%%%%%%%%%%%%%%%

\textbf{Notation.} We denote by  $\| \vx \| := \la \vx, \vx \ra^{1/2}$, $\vx \in \R^d$, the standard Euclidean norm for vectors,
and the spectral norm for symmetric matrices by $\| \mH \| := \max\{ \lambda_{\max}(\mH), - \lambda_{\min}(\mH) \}$, where $\mH = \mH^{\top} \in \R^{d \times d}$. For functions $h_1, h_2$, we write $h_1 = \cO(h_2)$ to denote $h_1 \leq C h_2$ for some universal numerical constant $C$, which does not depend
on the problem class parameters.

As is customary for second-order methods, we assume the objective function $f$ has $L$-Lipschitz Hessians, which is, for all $\vx, \vy \in \R^d$:
\begin{equation}\label{LipHess}
    \|\nabla^2 f(\vx) - \nabla^2f(\vy)\|\leq L\|\vx - \vy\|.
\end{equation}
We also assume that the objective is lower bounded from below, $f^\star:= \min_{\vx} f(\vx) > - \infty,$ and denote $F_0 = f(\vx_0) - f^\star$ for some initialization of the state $\vx_0$ that we fix throughout this paper without loss of generality.

The (Stochastic) Cubic Newton algorithm can be written as the following update rule of the current state~$\vx_t$:
%\begin{multline}
\beq \label{CNstep}
\ba{rcl}
    \vx_{t+1} &\in& \argmin_{\vy \in \R^d}
   \Bigl\{\, \Omega_{M,\vg,\mH}(\vy,\vx_t)\,\Bigr\},
\ea
\eeq
%\end{multline}
where $\Omega$
is the cubically regularized quadratic model of the objective, defined by
\beq \label{CubicModel}
\ba{cl}
&\Omega_{M, \vg, \mH}(\vy, \vx_t) 
:= \langle \vg_t, \vy - \vx_t \rangle \\[10pt]
& \;\;
   + \frac{1}{2} \langle \mH_t(\vy - \vx_t), \vy - \vx_t \rangle + \frac{M}{6}\|\vy - \vx_t\|^3,
\ea
\eeq
with regularization parameter
$M > 0$, which is typically greater than the Lipschitz constant $L$, whereas $\vg_t$ and $\mH_t$ are (stochastic) estimates of the gradient and Hessian at $\vx_t$.
Note that subproblem in~\eqref{CNstep}
is generally non-convex, when $\mH_t \not\succeq 0$,
and can have several isolated global minima.
However, it is possible to solve this problem efficiently
with the same arithmetic cost as for the classical Newton step.
In a non-degenerate case, one iteration of
the form~\eqref{CNstep}
can be represented in the standard form~\eqref{CubicNewtonStep},
regularizing the Hessian by identity matrix with a certain regularization coefficient $\gamma_t$,
which can be found by solving the corresponding univariate dual problem (see details in~\citep{conn2000trust,nesterov2006cubic,cartis2011adaptive1}).

\textbf{Examples of choices for $\vg_t$ and $\mH_t$:}
\begin{enumerate}
    \item Deterministic case: $\vg_t = \nabla f(\vx_t)\;, \mH_t = \nabla^2 f(\vx_t)$.
    \item Stochastic case, subsampled approaches: we sample batches $\cB_g,\cB_h$ of sizes $b_g,b_h$ respectively and set 
    \begin{equation}
        \label{subsamplingChoice}
    \ba{rcl}\vg_{t} &=& \frac{1}{b_g}\sum\limits_{i\in\cB_g }\nabla f_i(\vx_t), \\[10pt]
    \mH_{t} &=& \frac{1}{b_h}\sum\limits_{i\in\cB_h }\nabla^2 f_i(\vx_t).
    \ea
    \end{equation}
\end{enumerate}
In the deterministic case, Cubic Newton can be understood as replacing the function $f$ by its second-order Taylor approximation augmented by cubic regularization to control the region where this approximation is valid. When using stochastic samples, we further replace the exact Taylor model with the stochastic~\eqref{CubicModel}.

A point $\vx$ is denoted a \textit{$(\epsilon,c)$-approximate second-order local minimum} if it satisfies: 
$$
\ba{rcl}
\|\nabla f(\vx)\| & \leq & \epsilon
\quad \text{and} \quad
\lambda_{min}(\nabla^2 f(\vx)) \;\;\geq\;\; - c\sqrt{\epsilon},
\ea
$$ 
where $\epsilon, c > 0$ are given tolerance parameters.
Let us define the following accuracy measure (see \cite{nesterov2006cubic}):\vspace{-1mm}
$$
\ba{rcl}
\mu_c(\vx) := 
\max\Bigl(\|\nabla f(\vx)\|^{3/2}, \,
\frac{-\lambda_{min}(\nabla^2 f(\vx))^{3}}{c^{3/2}} \Bigr).
\ea
$$
Note that this definition implies
that if $\mu_c(\vx)\leq \epsilon^{3/2}$ then $\vx$ is an $(\epsilon,c)$-approximate 
\textit{local minimum}.

%%%%%%%%%%%%%%%%%%%%%%%%%%%%%%%%%%%%%%%%%%%%%
%%%%%%%%%%%%%%%%%%%%%%%%%%%%%%%%%%%%%%%%%%%%%
\subsection{Stochastic Cubic Newton}
\label{SectionStochCN}
%%%%%%%%%%%%%%%%%%%%%%%%%%%%%%%%%%%%%%%%%%%%%
%%%%%%%%%%%%%%%%%%%%%%%%%%%%%%%%%%%%%%%%%%%%%

For the sake of brevity, let us denote $\veps_t = \vg_t-\nabla f(\vx_t)$ (resp. $\mSigma_t = \mH_t - \nabla^2 f(\vx_t)$) the gradient (resp. Hessian) approximation error and $r_t = \|\vx_{t+1} - \vx_t\|$ the norm of the state update.

The main standard Lemma that we use for the analysis of one step of~\eqref{CNstep} is the following:
\begin{framedlemma} \label{lemma1} For any $M\geq L$, the progress made by the update rule \eqref{CNstep} satisfies:
    \begin{multline*}
        f(\vx_{t}) - f(\vx_{t+1}) \geq \cO(1)\Big( 
 \frac{1}{\sqrt{M}}\mu_M(\vx_{t+1}) \\
+ M r_{t}^3 - \frac{\|\veps_t\|^{3/2}}{\sqrt{M}} - \frac{\| \mSigma_t \|^3}{M^2}\Big)\;.
    \end{multline*}
\end{framedlemma}

In light of Lemma~\ref{lemma1}, we see that we need to control the $(3/2)-$th moment of the gradient noise $\|\veps_t\|^{3/2}$ and the third moment of the Hessian noise $\|\mSigma_t\|^3$, this is drastically different from the second moment usually encountered for first-order methods. The main complication that makes transferring methods from the first-order world to the second-order world is that we do not have an explicit expression of \eqref{CNstep}, making it very difficult to deal with the noise. In fact, Lemma~\ref{lemma1} gives the same bound if we assume the noise is deterministic (which is typically worse) or randomly centered, this suggests either a gap in theory or an intrinsic behavior of the method.
It remains an open question to establish stochastic lower bounds 
for second-order methods that enable distinguishing between these cases.

%Since the subject of lower bounds for second-order methods is still largely open, we, unfortunately, do not know which is which.

\textbf{Noise assumptions.} the usual noise assumptions made in the literature \citep{Subsampled1,Subsampled2,tripuraneni2017stochasticcubicregularizationfast,agafonov2024advancing} are the following:
\begin{framedassumption}[gradient noise]\label{GNoise}
    All samples~$\xi$ are sampled independently from the past, and there exists $\sigma_g<\infty$ such that$$
    \ba{cl}
    &\E \nabla f_\xi(\vx) \;\; = \;\; \nabla f(\vx)\;\;\;\text{and} \\[10pt]
    &\E\|\nabla f_\xi(\vx) - \nabla f(\vx)\|^2
    \;\; \leq \;\; \sigma_g^2\;.
    \ea
    $$
\end{framedassumption}
Note that by using Jensen's inequality, we have: $\E\|\nabla f_\xi(\vx) - \nabla f(\vx)\|^{3/2}\leq \sigma_g^{3/2}\;.$ By independence, the subsampled gradient in \eqref{subsamplingChoice} satisfies:
$$
\ba{rcl}
\E\|\veps_t\|^{3/2} &\leq& \frac{\sigma_g^{3/2}}{b_g^{3/4}},
\ea
$$
where $b_g$ is the batch size of the gradient samples.

Controlling the Hessian noise is much more complicated since we cannot use Jensen's inequality anymore, and it is easy to construct distributions with finite second-moment and infinite third moment. Thus we naturally need to make stronger assumptions.
\begin{framedassumption}[Hessian noise]\label{HessNoise}
    All samples $\xi$ are sampled independently from the past, and there exist positive constants $\sigma_h,\delta_h$ such that 
    $$
    \ba{cl}
    &\E \nabla^2 f_\xi(\vx) \;\; = \;\; \nabla^2 f(\vx), \\[10pt]
    &\E\|\nabla^2 f_\xi(\vx) - \nabla^2 f(\vx)\|^2
    \;\; \leq \;\; \sigma_h^2, 
    \;\; \text{and} \\[10pt]
    &\;\;\,\|\nabla^2 f_\xi(\vx) - \nabla^2 f(\vx)\|
    \;\; \leq \;\; \delta_h\;\;(a.s).
    \ea
    $$
    %\begin{multline*}
    %    \E \nabla^2 f_\xi(\vx) = \nabla^2 f(\vx),\E\|\nabla^2 f_\xi(\vx) - \nabla^2 f(\vx)\|^2\leq \sigma_h^2\\\text{and}\;\|\nabla^2 f_\xi(\vx) - \nabla^2 f(\vx)\|\leq \delta_h\;(a.s).
    %\end{multline*}
\end{framedassumption}
Assumption~\ref{HessNoise}, coupled with Theorem A.1 from \citep{chen2012maskedsamplecovarianceestimator} leads for example to the following bound on the subsampled Hessian in \eqref{subsamplingChoice}:
$$
\ba{rcl}
    \E\|\mSigma_t\|^3 
    &=& \cO\Bigl(\log(d)^{3/2}\frac{\sigma_h^3}{b_h^{3/2}} + \log(d)^3\frac{\delta_h^3}{b_h^3}\Bigr)\\[10pt]
    &\underset{b_h\gg 1}{\approx}& {\cO}\Bigl(\frac{\sigma_h^3}{b_h^{3/2}}\Bigr),
\ea
$$
where $b_h$ is the batch size of Hessian samples.

Under these assumptions, as shown by \citep{chayti2024unified}, using the estimates in \eqref{subsamplingChoice}, we get that for $T$ steps of \eqref{CNstep} and for $\vx_{out}^T$ sampled uniformly at random from the iterates $(\vx_t)_{1 \leq t\leq T}$:
\begin{equation}\label{PreviousRes}
   \E\mu_M(\vx_{out}^T) =\cO\left(\frac{\sqrt{M}F_0}{T} + \frac{\sigma_h^{3}}{M^{3/2} b_h^{3/2}} + \frac{\sigma_g^{3/2}}{b_g^{3/4}}\right).
\end{equation}
What is noteworthy about \eqref{PreviousRes} is that for $b_g = b_h = 1$, even with the best choice of $M$, \eqref{PreviousRes} can only guarantee convergence $\E\mu_M(\vx_{out}^T) \leq \varepsilon$ up to $\varepsilon \geq \cO(\sigma_g^{3/2})$.

\subsection{Momentum to the rescue}\label{Sec:MOM}

The main idea of this paper is to reuse past iterates to "forge'' a large batch. This concept is not new in the first-order optimization literature; The stochastic momentum or Polyak's momentum~\citep{POLYAK19641}, takes the form of an exponential moving average:
$$
\vg_t = (1 - \alpha_t) \vg_{t-1} + \alpha_t \nabla f_{\xi_t}(\vx_t)
\;\qquad\text{\color{blue}(HB)}
$$ 
for some $\alpha_t\in (0,1].$ 
This momentum is known to suffer from a high bias. To understand this, imagine $\vg_{t-1}= \nabla f(\vx_{t-1})$. Then, even in this perfect case, 
$$
\ba{rcl}
\!\!\!
\E\bigl[\vg_t - \nabla f(\vx_t)\bigr] & \!\!\!\! = \!\!\!\! &  
(1 - \alpha_t)(\nabla f(\vx_{t-1}) - \nabla f(\vx_{t}))\neq \vzero.
\ea
$$
Besides, its theoretical analysis usually involves assuming the gradient of 
$f$ is Lipschitz, an assumption that typically slows down
the Cubic Newton method,
and which is not required in our work
(see also the appendix, where we additionally discuss using the classical momentum
in the Cubic Newton).

%do not assume here (even if we did, the analysis of this momentum for Cubic Newton will not work which we will explain in the Appendix). 

At least three corrections were proposed for {\color{blue}(HB)}, the {\color{blue}MVR} (for Momentum based Variance Reduction) momentum \citep{cutkosky2020momentumbasedvariancereductionnonconvex}, the {\color{blue}IT} (stands for Implicit ``Gradient'' Transport) momentum \citep{arnold2019reducingvarianceonlineoptimization} and the second-order momentum ({\color{blue}SOM}) \cite{tran2022better}. %If we also count literature related to the Frank-Wolf algorithm, there is even a fourth momentum very close to ({\color{blue}SOM}) which we will denote  ({\color{blue}SOM'})\citep{zhang2019samplestochasticfrankwolfe}.

To describe these corrections,
let us consider that we observed a sequence of samples $(\xi_j)_{j \geq 0}$ and its corresponding estimates $(\pi_{\xi_j}(\vx_j))_{j \geq 0}$, which can be either sampled gradients or Hessians. Let $\alpha_t := 1/t$ and define the following quantity:
$\Pi_t(\vx) := \frac{1}{t}\sum_{j=0}^t \pi_{\xi_j}(\vx)\;. $

For illustration, let us assume that the functions $\pi_\xi$ are all affine,
$\pi_{\xi}(\vx) = \mA \vx + \vb_\xi$,
with the same matrix $\mA$ and different bias parts $\vb_\xi$ (e.g. $\pi_{\xi}$
are the gradients of a quadratic function).
%with the same gradient (
%$\exists A\; \forall\;\xi, \vx\;:\nabla \pi_\xi(\vx) = \mA$), 
Then we can easily show that
\begin{align*}
    \Pi_t(\vx_t)
    &= (1 - \alpha_t)\Big[\Pi_{t-1}(\vx_{t-1}) + \nabla \pi_{\xi_t}(\vx_t)(\vx_t - \vx_{t-1})\Big] \\&\quad+ \alpha_t \pi_{\xi_t}(\vx_t) \;\qquad\color{blue}\text{(SOM)}\\
    &= (1 - \alpha_t)\Big[\Pi_{t-1}(\vx_{t-1}) + \pi_{\xi_t}(\vx_t) - \pi_{\xi_t}(\vx_{t-1})\Big] \\&\quad+ \alpha_t \pi_{\xi_t}(\vx_t) \;\qquad\color{blue}\text{(MVR)}\\
    &= (1 - \alpha_t)\Pi_{t-1}(\vx_{t-1}) \\&\quad+ \alpha_t \pi_{\xi_t}\Big(\vx_t + \frac{1-\alpha_t}{\alpha_t}(\vx_t - \vx_{t-1})\Big).\;\qquad\color{blue}\text{(IT)}
\end{align*}
To have the corresponding expressions for the gradient (resp. Hessian) it suffices to take $\pi = \nabla f$ (resp. $\pi = \nabla^2 f$) and replace $\Pi_t(\vx_t)$ by $\vg_t$ (resp. $\mH_t$).

Although the three expressions we show are equivalent for affine $\pi$, they generally behave differently and need different assumptions to analyze their behavior.

Many combinations of momentum for the gradient and Hessian are possible but  we will only consider in the main paper, {\color{blue} (IT) } for the gradient and {\color{blue} (HB) } for the Hessian as this is the only combination that led to an improvement and most importantly to convergence for arbitrary batch sizes
for the stochastic Cubic Newton without needing any additional assumption. We consider other combinations in the appendix.

%%%%%%%%%%%%%%%%%%%%%%%%%%%%%%%%%%%%%%%%%%%%%%%%%%
%%%%%%%%%%%%%%%%%%%%%%%%%%%%%%%%%%%%%%%%%%%%%%%%%%
\section{Main theoretical analysis}
\label{SectionMainTheory}
%%%%%%%%%%%%%%%%%%%%%%%%%%%%%%%%%%%%%%%%%%%%%%%%%%
%%%%%%%%%%%%%%%%%%%%%%%%%%%%%%%%%%%%%%%%%%%%%%%%%%

Based on the previous discussion, we consider the following update rules for the gradient and Hessian estimates:
\begin{equation}\label{GradMom}
\ba{rcl}
    \vg_t &=& (1 - \alpha_t) \vg_{t-1} \\[5pt]
    && + \;\; \alpha_t \nabla f_{\xi_t}\bigl(\vx_t 
    + \frac{1-\alpha_t}{\alpha_t}(\vx_t - \vx_{t-1})\bigr),
\ea
\end{equation}
and 
\begin{equation}\label{HessMom}
    \mH_t = (1 - \beta_t) \mH_{t-1} + \beta_t \nabla^2 f_{\xi_t}(\vx_t ),
\end{equation}
where $\alpha_t,\beta_t\in(0,1]$.
The {\color{blue} (IT)} momentum in \eqref{GradMom}, is very similar to the Heavy Ball {\color{blue} (HB)} momentum and only differs in where the new gradient is evaluated. {\color{blue} (HB)} evaluates the new gradient at the current state $\vx_t$, while {\color{blue} (IT)} evaluates it at a ``transported'' location: $\vy_t = \vx_t + \frac{1-\alpha_t}{\alpha_t}(\vx_t - \vx_{t-1})$, thus these two moments have the same cost but {\color{blue} (IT)} is typically better.

Consequently, we consider the following algorithm:
\begin{algorithm}[h!]
	\caption{Cubic Newton with momentum}
	\label{Mainalg}
	\begin{algorithmic}[1]
		\Require $\vx_0 \in \R^d$, $M_t > 0$, $\alpha_t,\beta_t\in(0,1]$,  $T$. 
		\For{$t=0,\dotsc, T-1$}
        \State Sample $\xi^g_t, \xi^h_t$. $\qquad \text{ \# can be equal}$ 
        \If{$t = 0$}
        \State Compute $\vg_t = \nabla_{\xi^g_t} f(\vx_t), \mH_t = \nabla^2_{\xi^h_t} f(\vx_t)$.
        \Else
        \State Compute $\vg_t,\mH_t$ using \eqref{GradMom}, \eqref{HessMom}.
        \EndIf
		\State Take cubic step \\
        \quad\;\;$\vx_{t+1} \in \argmin_{\vy \in \R^d} \Omega_{M_t,\vg_t,\mH_t}(\vy,\vx_t)$
	\EndFor
\Return $\vx_{out}^T$ uniformly at random from $(\vx_i)_{0\leq i\leq T}$
	\end{algorithmic}
\end{algorithm}

%\nikita{Is it described somewhere how we could obtain $\vx_{out}^T$?}

To keep the analysis of Algorithm~\ref{Mainalg} simple, we consider our hyperparameters constant: the two ``momentum parameters'' $\alpha_t=\alpha,\beta_t=\beta$ and the regularization parameter $M_t=M$.

%\nikita{If we use only constant $\alpha$ and $\beta$ in our method, it could be even better to keep them constant everywhere in the text, to simplify things further (e.g., in (6) and in (HB)), as we do with parameter $M$. Then we might have a remark that potentially parameters can be adaptive, depending on $t$.}

\textbf{Note about noise Assumptions \ref{GNoise},\ref{HessNoise}.} We will assume we can bound the noise at the initial state~$\vx_0$ slightly differently, and denote $\sigma_{g,0}\leq \sigma_g, \sigma_{h,0}\leq \sigma_h$ its corresponding gradient and Hessian noise levels. We do this because, in some cases of our analysis (especially the convex case), the initial noise can give rise to slower terms in our convergence rates, such terms can be easily dealt with by making sure the initial gradient and Hessian estimates are of good quality (smaller variance), which can be done using large batches but only for the initial state. However, we stress that we do not need such a trick and can always guarantee convergence for any batch size even without it.

%%%%%%%%%%%%%%%%%%%%%%%%%%%%%%%%%%%%%%%%%%%%%%%%%%
%%%%%%%%%%%%%%%%%%%%%%%%%%%%%%%%%%%%%%%%%%%%%%%%%%
\subsection{Momentum bounds for Cubic Newton}
%%%%%%%%%%%%%%%%%%%%%%%%%%%%%%%%%%%%%%%%%%%%%%%%%%
%%%%%%%%%%%%%%%%%%%%%%%%%%%%%%%%%%%%%%%%%%%%%%%%%%

\textbf{Gradient momentum bound.} Based on the form of update~\eqref{GradMom}, we can write its error $\veps_t$ recursively as:
$$\veps_t = (1-\alpha)\veps_{t-1}+ (1-\alpha)\vz(\vx_{t-1},\vx_{t}) +\alpha\vz(\vy_t, \vx_t) + \alpha \vn_t\,,$$
where we define 
$$
\ba{rcl}\vn_t &:=& \nabla f_{\xi_t}(\vy_t) -\nabla f(\vy_t) \quad \text{and}\\[5pt] 
\vz(\va,\vb) &:=& \nabla f(\va) - \nabla f(\vb) - \nabla^2 f(\vb)(\va - \vb).
\ea
$$
Note that \eqref{LipHess} implies $\|\vz(\va,\vb)\| \leq \frac{L}{2}\|\va - \vb\|^2.$

\begin{framedlemma}\label{GradBound} For the gradient update defined in~\eqref{GradMom}, we show that the sequence $(\veps_t)_{t \geq 0}$ of gradient errors associated to it satisfies:
    \begin{multline*}
        \frac{1}{T}\sum_{t=0}^{T-1}\E\|\veps_t\|^{3/2}=\cO\Big(\overbrace{\alpha^{3/4}\sigma_g^{3/2}}^{\text{variance}} \\+ \underbrace{\frac{\sigma_{g,0}^{3/2}}{\alpha T}  +  \frac{(1-\alpha)^{3/2} L^{3/2}}{\alpha^{3}}\frac{1}{T}\sum_{t=0}^{T-1}\E r_{t}^{3}}_{bias}\Big)\,.
    \end{multline*}
\end{framedlemma}
Lemma~\ref{GradBound} shows that the average error of \eqref{GradMom} can be decomposed into two terms: \textbf{(variance)} term $\alpha^{3/4}\sigma_g^{3/2}$ and a \textbf{(bias)} $\frac{\sigma_{g,0}^{3/2}}{\alpha T}  +  \frac{(1-\alpha)^{3/2} L^{3/2}}{\alpha^{3}}\frac{1}{T}\sum_{t=0}^{T-1}\E r_{t}^{3}$. The \textbf{(variance)} term shows a smaller variance for a value of $\alpha$ that is small, while the \textbf{(bias)} term shows the opposite effect, this means that we should find a sweet spot to balance both terms. Note how the bias can also be decomposed into the \textbf{(initial bias)} $\frac{\sigma_{g,0}^{3/2}}{\alpha T}$ and the rest as \textbf{(iterate bias)} (iterates are different). The \textbf{initial bias} can be controlled using the initial gradient sample if needed. In contrast, the \textbf{iterate bias} can be controlled using the term $Mr_t^3$ that we left in Lemma~\ref{lemma1} (the main reason why other momentums do not work, in theory, is that they lead to an iterate bias with a power of $r_t$ different than $3$).

\textbf{Hessian momentum bound.}
Nearly the same story can be told about the Hessian momentum defined in~\eqref{HessMom}. We show that its error can be written recursively as:
$$
\ba{rcl}
\mSigma_t &=& (1 - \beta)\mSigma_{t-1}+ (1 - \beta)\mZ(\vx_{t-1},\vx_t) + \beta \mN_t,
\ea
$$
where 
$$
\ba{rcl}
\mN_t &:=& \nabla^2 f_{\xi_t}(\vy_t) -\nabla^2 f(\vy_t)
\quad \text{and} \\[5pt]
\mZ(\va,\vb) &:=& \nabla^2 f(\va) - \nabla^2 f(\vb).
\ea
$$
Again, \eqref{LipHess} implies that $\|\mZ(\va,\vb)\|\leq L \|\va - \vb\|$.

%\newpage

\begin{framedlemma}\label{HessBound}
For the Hessian update defined in \eqref{HessMom}, we show that the sequence $(\mSigma_t)_{t \geq 0}$ of Hessian errors associated to it satisfies:
\begin{multline*}
    \frac{1}{T}\sum_{t=0}^{T-1}\E\|\mSigma_t\|^3=\cO\biggl( \beta^{3/2} \Tilde{\sigma}_h^3(\beta) \\+ \frac{1}{\beta T}\sigma_{h,0}^3 +  \frac{(1 - \beta)L^3}{\beta^3}\frac{1}{T}\sum_{t=0}^{T-1} \E r_{t}^3\biggr)\,,
\end{multline*}
where 
$$
\ba{rcl}\Tilde{\sigma}_h^3(\beta)&:=&  \cO((\log(d)^{3/2}\sigma_h^3 + \log(d)^3\beta^{3/2}\delta_h^3 )\\[5pt]
&\underset{\beta\ll 1}{\approx}& {\cO}(\sigma_h^3).
\ea
$$
\end{framedlemma}
Lemma~\ref{HessBound} has the same structure as Lemma~\ref{GradBound} and the same remarks can be said about the bias and variance terms it contains.

To choose the momentum parameters $\alpha,\beta$, we proceed as follows:

\textbf{Controlling iterate bias terms. } As mentioned before, we use the term $M r_{t}^3$ in Lemma~\ref{lemma1} to compensate for the iterate bias terms we get in Lemmas~\ref{GradBound},~\ref{HessBound}.
We show we can get rid of this bias by choosing $\alpha,\beta \in(0,1]$ such that:
    \begin{equation}\label{Condition}
    \ba{rcl}
        \frac{4\sqrt{3} L^{3/2}}{\alpha^{3}} + \frac{73 }{M^{3/2}}\frac{9L^3}{\beta^3} - \frac{M^{3/2}}{72} & \leq & 0.
    \ea
    \end{equation}
To satisfy \eqref{Condition}, it suffices to take:
\begin{equation}\label{condition1}%ChoiceOfParams
\ba{rcl}
    \alpha &\geq& 10\sqrt{\frac{L}{M}} \quad\text{and}\quad 
    \beta \;\;\geq\;\; 46\frac{L}{M},
\ea
\end{equation}

\textbf{Controlling the initial bias terms. } To accomplish this, we further choose 
\begin{equation}\label{condition2}
    \ba{rcl}
        \alpha &\in& \underset{\alpha\in(0,1]}{\argmin}\;\;  \alpha^{3/4}\sigma_g^{3/2} + \frac{\sigma_{g,0}^{3/2}}{\alpha T} = \frac{a_g^{6/7}}{T^{4/7}}\\
    %\ea
    %\end{equation*} and \begin{equation*}
    %\ba{rcl}
    &\text{and}&\\
        \beta & \in &\underset{\beta\in(0,1]}{\argmin}\;\;   \beta^{3/2} \sigma_h^3 + \frac{1}{\beta T}\sigma_{h,0}^3 = \frac{a_h^{6/5}}{T^{2/5}},
    \ea
    \end{equation}
    
where, for simplicity, we denoted $a_g = \frac{\sigma_{g,0}}{\sigma_{g}}\leq 1$ and $a_h = \frac{\sigma_{h,0}}{\sigma_{h}}\leq 1$. 

All in all, based on \eqref{condition1} and \eqref{condition2}, we choose 
\begin{equation}\label{ChoiceOfParams}
\ba{rcl}
    \alpha = \max\Big(\frac{a_g^{6/7}}{T^{4/7}},10\sqrt{\frac{L}{M}}\Big)\, , \,
    \beta = \max\Big(\frac{a_h^{6/5}}{T^{2/5}},46\frac{L}{M}\Big)\,.
\ea
\end{equation}
To ensure $\alpha, \beta \leq 1$ it suffices to take $M\geq 100 L$.

We are ready to state our main convergence results.
\subsection{Non-convex rates}
In the general non-convex case, we show the following theorem:
\begin{framedtheorem}\label{NCVrate} For any $\alpha,\beta\in(0,1]$ and $M\geq L$ satisfying \eqref{Condition}, the iterates of Algorithm~\ref{Mainalg} are such that: 
    \begin{multline*}
        \E\mu_M(\vx_{out}) = \cO\biggl(\frac{\sqrt{M}F_0}{T} + \frac{ \sigma_{h,0}^3}{M^{3/2}\beta T} +  \frac{\sigma_{g,0}^{3/2}}{\alpha T} \\    + \frac{\beta^{3/2}\sigma_h^3}{M^{3/2}} + \alpha^{3/4}\sigma_g^{3/2} \biggr)\,.
    \end{multline*}
\end{framedtheorem}
As a sanity check, if we take $\alpha = \beta = 1$ (meaning we do not use momentum), we get the same result \eqref{PreviousRes} as past works \citep{chayti2024unified} (up to small $\cO(1/T)$ terms). The advantage of our approach is that we can pick $\alpha,\beta$ small that can satisfy the condition~\eqref{Condition} such as the choices in~\eqref{ChoiceOfParams}.

\begin{framedcorollary}\label{Corr}
For $\alpha,\beta$ defined in \eqref{ChoiceOfParams}, and any $M\geq 100 L$, we have:
\begin{multline*}
        \E\mu_M(\vx_{out}) = \cO\biggl(\frac{\sqrt{M}F_0}{T}  + \frac{L^{3/2}\sigma_h^3}{M^{3}} + \frac{L^{3/8}}{M^{3/8}}\sigma_g^{3/2}\\ + \frac{a_{g}^{9/14}\sigma_{g}^{3/2}}{T^{3/7}} + \frac{a_h^{9/5}\sigma_{h}^{3}}{M^{3/2} T^{3/5}} \biggr).
    \end{multline*}
\end{framedcorollary}
%Note that the quantity $\Tilde{F}_0 $ in Corrolary~\ref{Corr} is slightly different from the usual constant $F_0$, but their difference is bounded by a constant that depends on the initial gradient and Hessian noise, which means that we can make it as small as possible by for example assuming the first batch is proportional to $T$ in size (e.g. $T/100$); however, we do not have to do this here. 

The last two terms represent the effect of the initial gradient and Hessian noise on convergence and we can see that they decay with $T$.

By choosing $M$ such that it minimizes the right-hand side in Corollary~\ref{Corr}, we get:
\begin{multline}\label{NCVrate}
    \E\mu_M(\vx_{out})=\cO\biggl(\frac{\sqrt{L}F_0}{T}  \\ + \frac{L^{3/14}F_0^{6/7}\sigma_h^{3/7}}{T^{6/7}}+ \frac{L^{3/14}F_0^{3/7}\sigma_g^{6/7}}{T^{3/7}}\\
    + \frac{a_{g}^{9/14}\sigma_{g}^{3/2}}{T^{3/7}} + \frac{a_h^{9/20}\sigma_{h}^{3/4}F_0^{1/2}}{ T^{9/10}}\biggr)\,.
\end{multline}

The best old rate in $T$ is $\cO(\frac{\sqrt{L}F_0}{T}  + \big(\frac{F_0\sigma_h}{T}\big)^{3/4}+ \sigma_g^{3/2})$.
Note that we improve both the $\sigma_h$ term from $\cO(T^{-3/4})$ to $\cO(T^{-6/7})$ and, most importantly, we change the constant term $\cO(\sigma_g^{3/2})$ to a decreasing term $\cO(T^{-3/7})$, thus guaranteeing convergence for any batch size.

%%%%%%%%%%%%%%%%%%%%%%%%%%%%%%%%%%%%%%%%%%%%%%%%%%
%%%%%%%%%%%%%%%%%%%%%%%%%%%%%%%%%%%%%%%%%%%%%%%%%%
\subsection{Extension to the convex case}
\label{SectionConvex}
%%%%%%%%%%%%%%%%%%%%%%%%%%%%%%%%%%%%%%%%%%%%%%%%%%
%%%%%%%%%%%%%%%%%%%%%%%%%%%%%%%%%%%%%%%%%%%%%%%%%%

Let us consider the case when the objective function $f$ is convex\footnote{Note that we cover, in the appendix, the more general composite problems: $\min_{\vx\in\R^d} F(\vx):= \E_{\xi\sim \mathcal{P}}\big[f_\xi(\vx)\big] + \psi(\vx)$,
where $\psi$ is a simple possibly nonsmooth convex function. \
}. We assume that all iterates $(\vx_t)_{t \geq 0}$
are bounded within a ball of radius $D$: $\| \vx_t - \vx^{\star} \| \leq D$,
where $\vx^{\star}$ is any of the global minimizers
that we assume to exist.
This assumption is customary and can be easily achieved, for example, by introducing an additional ball constraint to our problem.
These assumptions provide us with the following inequality:
\vspace*{-5pt}
\begin{equation} \label{ConvexityBound}
    f(\vx) - f^\star \leq \la \nabla f(\vx), \vx - \vx^{\star} \ra \leq  D \|\nabla f(\vx)\|\,,
\end{equation}
which means that convergence to a stationary point implies convergence to a global minimum.
Under inequality~\eqref{ConvexityBound}, previous results \citep{TensorCN,chayti2024unified} 
give the following rate for the Stochastic Cubic Newton:
%(without acceleration) give the following rate:
\begin{equation}\label{ConvexStoch}
    \E[f(\vx^T_{out})] - f^\star =  \cO\Big(\frac{LD^3}{T^2} + \frac{\sigma_h D^2}{T} + \sigma_g D\Big).
\end{equation}
We observe the same problem as in the non-convex case:
the method converges only to a neighborhood of the solution, and it cannot achieve accuracy $\varepsilon \leq \sigma_g D$.%, for small batches.
%For the sake of simplicity, we will assume the initial noise constants $\sigma_{g,0}, \sigma_{h,0}$ are small enough that we can ignore them (the full, complicated statement including these terms can be found in the Appendix). 
%We can show the following Theorem:
\begin{framedtheorem}\label{ConvTH} For $\alpha,\beta$ defined in \eqref{ChoiceOfParams} and $M\geq 100 L$, when the objective $f$ is convex, the iterates of Algorithm~\ref{Mainalg} satisfy:
\begin{multline*}
\min_{0\leq t\leq T-1} F_t=\cO\Big( \frac{M D^3}{T^2} + \frac{L D \Tilde{\sigma}_h^2}{M^2} + \frac{L^{1/4}}{M^{1/4}}D\sigma_g\\ + F_0 e^{-T/3} + \frac{a_g^{3/7} D \sigma_g}{T^{2/7}} + \frac{a_h^{6/5} D \sigma_h^2}{MT^{2/5}} \Big),
\end{multline*}
where we define $F_t = \E[f(\vx_t)] - f^\star$.

Choosing $M$ such that it minimizes the right-hand side, we get:
\begin{multline*}
\min_{0\leq t\leq T-1} F_t\\=\cO\Big( \frac{L D^3}{T^2} + \frac{L^{1/3}D^{7/3}\Tilde{\sigma}_h^{2/3}}{T^{4/3}} + \frac{L^{1/5}D^{7/5}\sigma_g^{4/5}}{T^{2/5}}\\ + F_0 e^{-T/3} + \frac{a_h^{3/5} D^{2}\sigma_h}{T^{6/5}} + \frac{a_g^{3/7} D \sigma_g}{T^{2/7}}  \Big).
\end{multline*}
\end{framedtheorem}
The last three terms in Theorem~\ref{ConvTH} correspond to the effect of the initialization and the initial gradient and Hessian noises. The last two terms can be made even smaller by choosing a larger initial batch (i.e. reducing $a_g,a_h$).

Again, we improve both Hessian and gradient rates and can converge for any precision $\varepsilon>0$ for any batch size.

By ignoring the last three terms, we get the following iteration complexity:
$\cO\Big((\frac{LD^3}{\varepsilon})^{1/3} + \frac{L^{1/4}D^{7/4}\Tilde{\sigma}_h^{1/2}}{\varepsilon^{3/4}} + \frac{L^{1/2}D^{7/2}\sigma_g^{2}}{\varepsilon^{5/2}}\Big)$.

We note that in the convex case, it is possible to reach
better rates for the Cubic Newton using acceleration~\cite{nesterov2008accelerating}.
For example, \citep{agafonov2024advancing} have a rate 
for their stochastic second-order method:
$\cO\Big((\frac{LD^3}{\varepsilon})^{1/3} + \frac{D\Tilde{\sigma}_h^{1/2}}{\varepsilon^{1/2}} + \frac{D\sigma_g}{\varepsilon^{2}}\Big)$. We believe our momentum can be coupled with acceleration as well, and has the potential to improve this result further.

\section{Experiments}

\begin{figure}
    \vspace*{-12pt}
    \centering
    \includegraphics[width=0.5\linewidth]{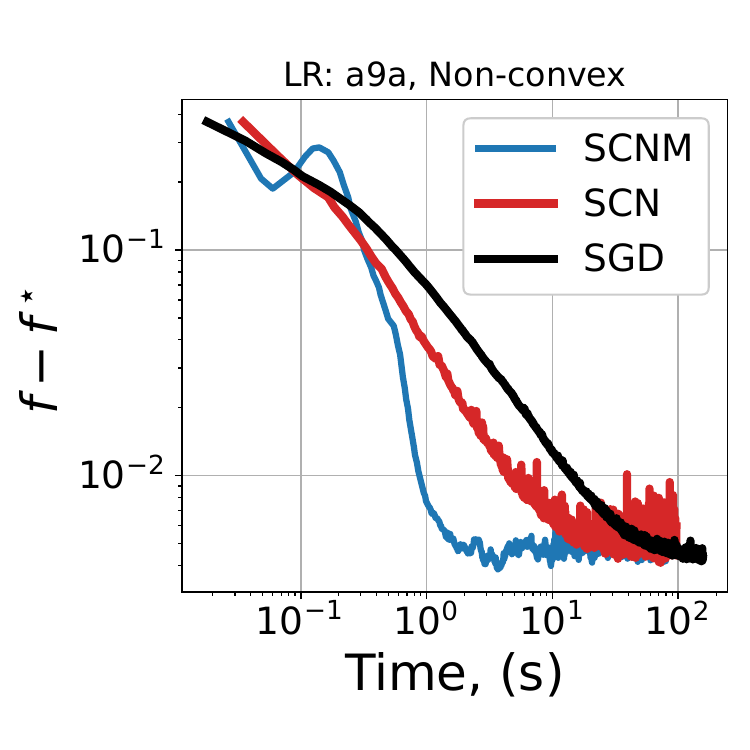}
    \caption{\small A small batch size $b_g=512, b_h=216$ is used. We see SCNM, our Algorithm with momentum, outperforms its counterpart without momentum, it also outperforms SGD.}
    \label{fig:1}
\end{figure}

\begin{figure}
    \centering
    \includegraphics[width=0.5\linewidth]{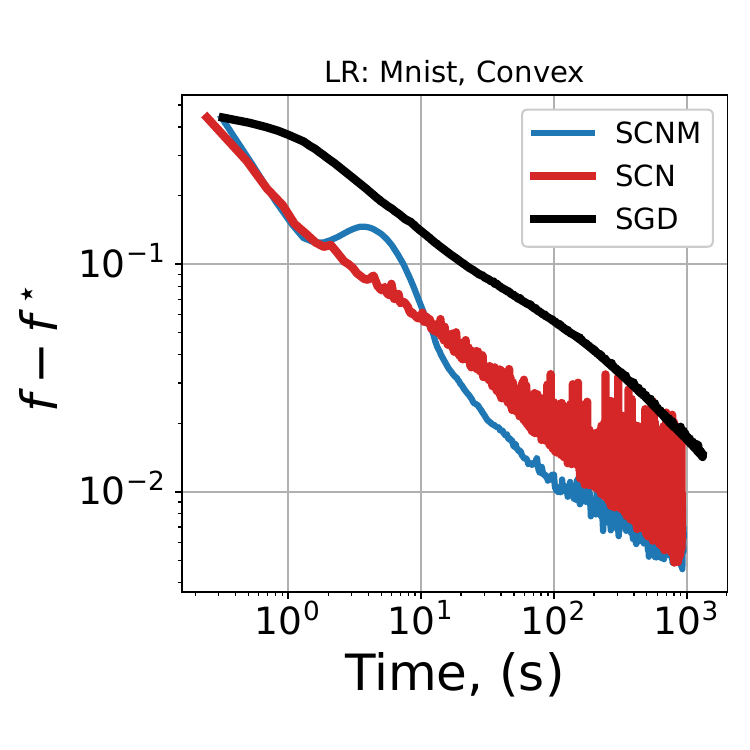}
    \caption{\small Logistic regression on the MNIST dataset. A small batch size $b_g=512, b_h=216$ is used. SCNM, our Algorithm with momentum, again outperforms its counterpart without momentum and SGD.}
    \label{fig:2}
\end{figure}

We consider a logistic regression with a non-convex regularization ($Reg(\vx) = \sum_{i=1}^d \frac{\vx_i^2}{1+\vx_i^2}$), following~\cite{Subsampled1}. We use the A9A dataset from the LibSVM library \citep{chang2011libsvm}. Results are shown in Figure~\ref{fig:1}.
We also consider the MNIST dataset \citep{mnist} and show the results in Figure~\ref{fig:2}.
Both experiments depicted in Figures~\ref{fig:1},\ref{fig:2} use momentum parameter values $\alpha = 0.1$ and $\beta=0.01$.
We see that using momentum in the Stochastic Cubic Newton indeed improves its convergence rate and visibly decreases the variance, as predicted by our theory.

\vspace*{-5pt}

\section{Discussion and Limitations}

%\textbf{Potential applications beyond this work.}
%We showed in this work how to use a specific type of momentum to improve the rate of convergence of the Stochastic Cubic Newton method and most importantly make it converge for any choice of the batch sizes. We believe this idea can, potentially, be combined with other ideas previously used in the literature like acceleration and quadratic regularization, which were shown in \citep{agafonov2024advancing} to help improve the gradient noise term. More generally, momentum can be combined with other second-order methods like the Levenberg-Marquardt regularization \citep{mishchenko2023regularized} for convex functions, which still lacks a stochastic version.  

%\nikita{Maybe we should keep it for ourselves and not advertise broadly for now? I believe acceleration could be obtained relatively easy}
\textbf{Practical implementation.} We did not provide a specific method for solving the cubic sub-problem in Algorithm~\ref{Mainalg}. A naive implementation would require access to the entire Hessian matrix, making it impractical for high-dimensional problems. However, we can solve this problem using only Hessian-vector products, which are typically more efficient for various problem classes, such as neural networks with backpropagation. It is important to note that incorporating Hessian momentum may complicate this process. Therefore, we leave the exploration of how to effectively combine Hessian-vector products and momentum for future work. 

\textbf{Other momentum choices.} As shown in Subsection~\ref{Sec:MOM}, there are other types of momentum that we tried in the supplementary material. We were able to improve the gradient noise term from $\cO(\frac{1}{T^{3/7}})$ in \eqref{NCVrate} to $\cO(\frac{1}{T^{9/20}})$ using the {\color{blue}MVR} momentum from \citep{cutkosky2020momentumbasedvariancereductionnonconvex}; this improvement comes at the cost of needing to assume the gradient of the objective is $L_g$-Lipschitz, but this constant can be made as big as we want by adding quadratic functions to our objective. Even when ignoring this problem with the constant $L_g$, the corresponding iteration complexity is only $\cO(\frac{1}{\varepsilon^{30/9}})$, which is worse than $\cO(\frac{1}{\varepsilon^{3}})$ obtained using the corresponding momentum in first-order methods \citep{cutkosky2020momentumbasedvariancereductionnonconvex}. One of the reasons we think this is the case is that we did not have corresponding powers of $r_t$ in the iterate bias that let us compensate for the bias using the term $Mr_t^3$ in Lemma~\ref{lemma1}. We note that this exact problem is present for first-order methods (although it is less complicated since we know $r_t$ explicitly in this case) and the fix was to use normalized methods such as in \citep{cutkosky2020momentumimprovesnormalizedsgd} to make the size of the step constant. Such a technique is still absent for second-order methods, and were it to be discovered, it would mean the possibility of improving the other types of momentum we discussed in this work.

\textbf{Adaptive versions.} The main drawback of our method is that it introduces an additional two hyperparameters $\alpha, \beta$; while we can choose typical values as we did in the experiments, developing strategies to adapt to these new hyperparameters is important. We see our work as a first step in showcasing the potential benefits of using momentum in second-order methods
and keep these questions for future research.

\section{Conclusion}
In this work, we study Stochastic Cubic Newton (SCN) with a specific type of momentum for non-convex problems and show improved rates compared to the plain version without momentum. Thus, we show, for the first time, the convergence of SCN for arbitrary batch sizes.

\section*{Acknowledgments}

This work was supported by the Swiss State Secretariat for
Education, Research and Innovation (SERI) under contract
number 22.00133.

\bibliography{main}

\appendix
\onecolumn
\section{Additional results using other momentums}
We have also considered other combinations of momentum for the gradient estimation in the Cubic Newton method. We found out that they all guarantee convergence to arbitrary accuracy for any batch size. However, they make the Hessian term in the convergence rate slower, and only in the case of the MVR momentum we are able to show an improved rate for the gradient term.

We also note that to analyze the HB and MVR momentums, we needed to make additional assumptions on the smoothness of the objective function $f$. Specifically, for HB, we assumed that the gradients of $f$ are $L_g$-Lipschitz, and for MVR, we needed to assume that the gradients of the sample functions $(f_{\xi})$ are all $L_g$-Lipschitz.

\begin{align*}
    \text{(SOM) Section~\ref{sec:SOM}:}\quad &\vg_t = (1 - \alpha) \Big(\vg_{t-1} + \mH_t(\vx_t - \vx_{t-1})\Big) + \alpha \nabla f_{\xi_t}(\vx_t ),\\
    \text{(HB) Section~\ref{sec:HB}:}\quad &\vg_t = (1 - \alpha) \vg_{t-1} + \alpha \nabla f_{\xi_t}(\vx_t ),\\
    \text{(MVR) Section~\ref{sec:MVR}:}\quad &\vg_t = (1 - \alpha) \Big(\vg_{t-1} + \nabla_{\xi_t} f(\vx_t) - \nabla_{\xi_t} f(\vx_{t-1})\Big) + \alpha \nabla_{\xi_t} f(\vx_t),
\end{align*}
and for all these gradient estimators, we consider the HB momentum for the Hessian:
$$\mH_t = (1 - \beta) \mH_{t-1} + \beta \nabla^2 f_{\xi_t}(\vx_t), \;\text{for}.$$

To give a brief summary of the new results, let us assume for simplicity that $\sigma_{g,0} = 0$. Then we get the following rates after $T$ iterations of our algorithm using these different types of momentum:

$$
\E\mu_M(\vx_{out}) = \left\{
    \begin{array}{ll}
        \cO\Big( \frac{\sqrt{L}F_0}{T} + \Big(\frac{F_0\sigma_h}{T}\Big)^{3/4} + \frac{\sigma_g^{9/16}\sigma_h^{3/8} F_0^{3/8}}{T^{3/8}} + \frac{L^{3/10}\sigma_g^{3/5} F_0^{3/5}}{T^{3/5}}  \Big), & \mbox{(SOM) Section~\ref{sec:SOM} } \\
        \cO\Big(\frac{\sqrt{L}F_0}{T}+ \Big(\frac{F_0\Tilde{\sigma}_h}{T}\Big)^{3/4} + \Big(\frac{L_g F_0\sigma_g^2}{T}\Big)^{3/8} \Big) , & \mbox{(HB) Section~\ref{sec:HB}}\\
        \cO\Big(\frac{\sqrt{L}F_0}{T} + \Big(\frac{F_0\Tilde{\sigma}_h}{T}\Big)^{3/4} + \Big(\frac{L_g F_0}{T}\Big)^{9/20}\sigma_g^{6/10}\Big) , & \mbox{(MVR) Section~\ref{sec:MVR}}
    \end{array}
\right.
$$
where $\Tilde{\sigma}_h = \cO(L_g + \sigma_h)$.

In particular, we see an improved dependence on $T$ for the ``gradient'' term of the (MVR) momentum which is better than the $\frac{1}{T^{3/7}}$ that we got in the main paper with the (IT) momentum. However, this improved dependence comes at the cost of needing an additional assumption and a dependence on $L_g$ instead of $L$ (note that $L_g$ can be made arbitrarily big by adding a quadratic to our objective function, while $L$ will be unchanged).

\section{Preliminaries}
We consider the general problem 
$$\min\limits_{\vx\in\R^d} f(\vx) := \E_\xi[f_\xi(\vx)]\;,$$

\noindent where $f_\xi$ is twice differentiable and $f$ has $L$-Lipschitz Hessian i.e.:\begin{equation}\label{AppLipHess}
    \| \nabla^2 f(\vx) - \nabla^2 f(\vy) \|  \leq  L \|\vx - \vy\|, \qquad \forall \vx, \vy \in \R^d .
\end{equation}
As a direct consequence of \eqref{AppLipHess} (see \citep{nesterov2006cubic,nesterov2018lectures}) we have for all $\vx,\vy\in\R^d$: \begin{equation}\label{LipHessGrad}
    \| \nabla f(\vy) - \nabla f(\vx) - \nabla^2 f(\vx)(\vy - \vx)\|  \leq  \frac{L}{2}\|\vx - \vy\|^2,
\end{equation} 
\begin{equation}\label{LipHessFunc}
    | f(\vy) - f(\vx) - \la \nabla f(\vx), \vy - \vx \ra - \frac{1}{2} \la \nabla^2 f(\vx)(\vy - \vx), \vy - \vx \ra |
 \leq  \frac{L}{6}\|\vy - \vx\|^3.
\end{equation}
For a given vector $\vg$ and matrix $\mH$, we define $\vx^{+}$ as 
\begin{equation}\label{xx+}
\ba{rcl}
\!\!\!\!\!\!
    \vx^+ & \in & \argmin_{\vy \in \R^d}
    \Bigl\{\, \Omega_{M,\vg,\mH}(\vy,\vx) := 
\langle \vg, \vy - \vx \rangle + \frac{1}{2} \langle \mH(\vy - \vx), \vy - \vx \rangle
+ \frac{M}{6}\|\vy - \vx\|^3 \,
\Bigr\}.
\ea
\end{equation}

We define the quantity $$\mu_M(\vx) = \max(\|\nabla f(\vx)\|^{3/2}, \frac{-\lambda_{min}(\nabla^2 f(\vx))^{3}}{M^{3/2}}),$$
where $\lambda_{min}(\mA)$ denotes the smallest eigenvalue of the matrix $\mA$.

\cite{chayti2024unified} show the following theorem (see Appendix C, Theorem C.1):
\begin{framedtheorem} \label{ThOneStep}
    \label{inexactCN}
	For any $\vx \in \R^d$,
	let $\vx^+$ be defined by \eqref{xx+}. Then, for $M\geq L$ we have:
\begin{equation*}
\ba{rcl}
    \label{MainInequality}
f(\vx) - f(\vx^+) & \geq & 
 \frac{1}{1008\sqrt{M}}\mu_M(\vx^+) 
+ \frac{M \| \vx - \vx^+\|^3 }{72} - \frac{4 \|\nabla f(\vx) - \vg\|^{3/2}}{\sqrt{M}} - \frac{73 \| \nabla^2 f(\vx) - \mH \|^3}{M^2}.
\ea
\end{equation*}
\end{framedtheorem}

We will also need the following Lemma to deal with averages under cubic powers:
\begin{framedlemma}
    \label{lyapMatr}
  Suppose that $q \geq 2$, $p\geq 2$, and fix $r\geq \max{(q,2\log(p))}$. Consider (weakly) independent random
self-adjoint matrices $Y_1, \cdots, Y_N$ with dimension $p \times p$, $\E[Y_i] = 0$. It holds that: 
$$
\big[\E[\|\sum_{i=1}^N Y_i\|^q]\big]^{1/q} \leq 2\sqrt{er} \|\big(\sum_{i=1}^N \E[Y_i^2]\big)^{1/2}\| + 4 er \E[\max_i \|Y_i\|^q]^{1/q}.
$$
\end{framedlemma}
The proof of this Lemma follows exactly the same lines as Lemma 32 from \citep{SVRC} and Theorem A.1 from \citep{chen2012maskedsamplecovarianceestimator}, and noticing that for the symmetrization trick to work, we only need the sequence to be (weakly) independent.

\section{Proofs of results of the main paper}

\subsection{Different forms of momentum}

In the paper, we introduce the quantity $$\Pi_t(\vx) := \frac{1}{t}\sum_{j=0}^t \pi_{\xi_j}(\vx)\;. $$
and define $\alpha_t = 1/t$. 

Our aim is, at each step $t$, to approximate $\Pi_t(\vx_t)$.

We also assume for illustration purposes that the functions $\pi_\xi$ are all linear with the same gradient $\mA$.

Then we have:
\begin{align*}
    \Pi_t(\vx_t)
    &= (1 - \alpha_t)\Pi_{t-1}(\vx_{t})+ \alpha_t \pi_{\xi_t}(\vx_t) \\
    &= (1 - \alpha_t)\Big[\Pi_{t-1}(\vx_{t-1}) + \mA (\vx_t - \vx_{t-1})\Big] + \alpha_t \pi_{\xi_t}(\vx_t) \\
    &= (1 - \alpha_t)\Big[\Pi_{t-1}(\vx_{t-1}) + \nabla \pi_{\xi_t}(\vx_t)(\vx_t - \vx_{t-1})\Big] \quad+ \alpha_t \pi_{\xi_t}(\vx_t) \;\qquad\color{blue}\text{(SOM)}\\
    &= (1 - \alpha_t)\Big[\Pi_{t-1}(\vx_{t-1}) + \pi_{\xi_t}(\vx_t) - \pi_{\xi_t}(\vx_{t-1})\Big] \quad+ \alpha_t \pi_{\xi_t}(\vx_t) \;\qquad\color{blue}\text{(MVR)}\\
    &= (1 - \alpha_t)\Pi_{t-1}(\vx_{t-1}) + \alpha_t \Big(\pi_{\xi_t}(\vx_t) + \frac{1-\alpha_t}{\alpha_t}\nabla \pi_{\xi_t}(\vx_t)(\vx_t - \vx_{t-1})\Big)\\
    &= (1 - \alpha_t)\Pi_{t-1}(\vx_{t-1}) + \alpha_t \pi_{\xi_t}\Big(\vx_t + \frac{1-\alpha_t}{\alpha_t}(\vx_t - \vx_{t-1})\Big).\;\qquad\color{blue}\text{(IT)}
\end{align*}

\subsection{Momentum bounds}
The gradient and Hessian updates we used in the main paper have the following form: 

\begin{equation}\label{Gradmom}
    \vg_t = (1 - \alpha) \vg_{t-1} + \alpha \nabla f_{\xi_t}\Big(\vx_t + \frac{1-\alpha}{\alpha}(\vx_t - \vx_{t-1})\Big), \;\text{for}\;t\geq 1\;\text{and}\; \vg_0 = \nabla f_{\xi_0}(\vx_0)
\end{equation}

and 
\begin{equation} \label{Hmom}
    \mH_t = (1 - \beta) \mH_{t-1} + \beta \nabla^2 f_{\xi_t}(\vx_t), \;\text{for}\;t\geq 1\;\text{and}\; \mH_0 = \nabla^2 f_{\xi_0}(\vx_0)\;,
\end{equation}
for $\alpha,\beta\in(0,1]$.

We define $r_t = \|\vx_t - \vx_{t-1}\|$.

\textbf{Bounding the gradient momentum error.}
Let $\veps_t := \vg_t - \nabla f(\vx_t)$ and $\vz(\va,\vb) := \nabla f(\va) - \nabla f(\vb) - \nabla^2 f(\vb)(\va - \vb)$.

Note that \eqref{LipHessGrad} implies $\|\vz(\va,\vb)\| \leq \frac{L}{2}\|\va - \vb\|^2.$

For ease of notation, we will also define $\vtheta_t := \vx_t + \frac{1-\alpha}{\alpha}(\vx_t - \vx_{t-1})$, $\vn_t = \nabla f_{\xi_t}(\vtheta_t)  - \nabla f(\vtheta_t)$ and $\vn_0 = \nabla f_{\xi_0}(\vx_0)  - \nabla f(\vx_0)$.
\begin{framedlemma}
For $t\geq 1$:
    \begin{align*}
    \frac{1}{T}\sum_{t=0}^{T-1}\E\|\veps_t\|^{3/2}\leq \frac{3}{\alpha T} \sigma_{g,0}^{3/2} + \sqrt{3} \alpha^{3/4}\sigma_g^{3/2} + \frac{\sqrt{3}(1-\alpha)^{3/2} L^{3/2}}{\alpha^{3}}\frac{1}{T}\sum_{t=0}^{T-1}\E r_{t}^{3}
\end{align*}
\end{framedlemma}
\begin{proof}

For $t\geq 1$, we have:
\begin{align*}
    \veps_t &= (1-\alpha) \vg_{t-1} + \alpha \nabla f_{\xi_t}(\vtheta_t) - \nabla f(\vx_t) \\
    &= (1-\alpha) \Big(\vg_{t-1} - \nabla f(\vx_t)\Big) + \alpha \Big(\nabla f_{\xi_t}(\vtheta_t) - \nabla f(\vx_t)\Big) \\
     &= (1-\alpha) \Big(\veps_{t-1}+ \nabla f(\vx_{t-1}) - \nabla f(\vx_{t})\Big) + \alpha \Big(\vn_t + \nabla f(\vtheta_t) - \nabla f(\vx_t)\Big) \\
     &= (1-\alpha)\Big(\veps_{t-1}+ \vz(\vx_{t-1},\vx_{t})\Big) + \alpha \Big(\vn_t + \nabla f(\vtheta_t) - \nabla f(\vx_t) + \frac{1-\alpha}{\alpha}\nabla^2 f(\vx_t)(\vx_{t-1} - \vx_t)\Big) \\
     &= (1-\alpha)\Big(\veps_{t-1}+ \vz(\vx_{t-1},\vx_{t})\Big) + \alpha \Big(\vn_t + \nabla f(\vtheta_t) - \nabla f(\vx_t) - \nabla^2 f(\vx_t)(\vtheta_t - \vx_t)\Big) \\
     &= (1-\alpha)\Big(\veps_{t-1}+ \vz(\vx_{t-1},\vx_{t})\Big) + \alpha \Big(\vn_t +\vz(\vtheta_t, \vx_t)\Big) \\
     &= (1-\alpha)\veps_{t-1}+ (1-\alpha)\vz(\vx_{t-1},\vx_{t}) +\alpha\vz(\vtheta_t, \vx_t) + \alpha \vn_t. \\
\end{align*}

Thus
$$\veps_t = (1 - \alpha)^t\veps_0 + \sum_{i=0}^{t-1} (1-\alpha)^i\Big((1-\alpha)\vz(\vx_{t-i-1},\vx_{t-i}) +\alpha\vz(\vtheta_{t-i}, \vx_{t-i})\Big) + \alpha \sum_{i=0}^{t-1} (1-\alpha)^i\vn_{t-i}\;,\; t\geq 1,$$

and 
\begin{align*}
    \|\veps_t\|&\leq (1 - \alpha)^t\|\veps_0\| + \sum_{i=0}^{t-1} (1-\alpha)^i\Big((1-\alpha)\|\vz(\vx_{t-i-1},\vx_{t-i})\| +\alpha\|\vz(\vtheta_{t-i}, \vx_{t-i})\|\Big) + \alpha\|\sum_{i=0}^{t-1} (1-\alpha)^i\vn_{t-i}\|\\
    &\leq (1 - \alpha)^t\|\veps_0\| + \frac{L}{2}\sum_{i=0}^{t-1} (1-\alpha)^i\Big((1-\alpha)\|\vx_{t-i-1}-\vx_{t-i}\|^2 +\alpha\|\vtheta_{t-i}-\vx_{t-i}\|^2\Big) + \alpha \|\sum_{i=0}^{t-1} (1-\alpha)^i\vn_{t-i}\|\\
    &\leq (1 - \alpha)^t\|\veps_0\| + \frac{L}{2}\sum_{i=0}^{t-1} (1-\alpha)^i\Big((1-\alpha)r_{t-i}^2 +\frac{(1-\alpha)^2}{\alpha}r_{t-i}^2\Big) + \alpha\|\sum_{i=0}^{t-1} (1-\alpha)^i\vn_{t-i}\|\\
    &\leq (1 - \alpha)^t\|\veps_0\| + \frac{(1-\alpha) L}{\alpha}\sum_{i=0}^{t-1} (1-\alpha)^ir_{t-i}^2 + \alpha\|\sum_{i=0}^{t-1} (1-\alpha)^i\vn_{t-i}\|.\\
\end{align*}

Now we raise this to the power $3/2$ and using the convexity of $x\mapsto x^{3/2}$, we get:
$$\|\veps_t\|^{3/2}\leq \sqrt{3}\Bigg((1 - \alpha)^{3t/2}\|\veps_0\|^{3/2} + \frac{(1-\alpha)^{3/2} L^{3/2}}{\alpha^{3/2}}\Big(\sum_{i=0}^{t-1} (1-\alpha)^i r_{t-i}^2\Big)^{3/2} + \alpha^{3/2}\|\sum_{i=0}^{t-1} (1-\alpha)^i\vn_{t-i}\|^{3/2}\Bigg).$$
Again by the convexity of $x\mapsto x^{3/2}$, we have:
\begin{align*}
    \Big(\sum_{i=0}^{t-1} (1-\alpha)^i r_{t-i}^2\Big)^{3/2}&\leq \Big(\sum_{i=0}^{t-1} (1-\alpha)^i \Big)^{1/2} \sum_{i=0}^{t-1} (1-\alpha)^i r_{t-i}^{3}\\
    &\leq (\frac{1}{\alpha})^{1/2} \sum_{i=0}^{t-1} (1-\alpha)^i r_{t-i}^{3}.
\end{align*}

Thus $$\|\veps_t\|^{3/2}\leq \sqrt{3}\Bigg((1 - \alpha)^{3t/2}\|\veps_0\|^{3/2} + \frac{(1-\alpha)^{3/2} L^{3/2}}{\alpha^{2}}\sum_{i=0}^{t-1} (1-\alpha)^i r_{t-i}^{3} + \alpha^{3/2}\|\sum_{i=0}^{t-1} (1-\alpha)^i\vn_{t-i}\|^{3/2}\Bigg).$$
Now we take the expectation over $(\xi_j)_{j\leq t}$, we get:
$$\E\|\veps_t\|^{3/2}\leq \sqrt{3}\Bigg((1 - \alpha)^{3t/2}\E\|\vn_0\|^{3/2} + \frac{(1-\alpha)^{3/2} L^{3/2}}{\alpha^{2}}\sum_{i=0}^{t-1} (1-\alpha)^i \E r_{t-i}^{3} + \alpha^{3/2}\E\|\sum_{i=0}^{t-1} (1-\alpha)^i\vn_{t-i}\|^{3/2}\Bigg).$$

Assuming $\E \nabla f_\xi(\vx_j) = \nabla f_\xi(\vx_j)$ and $\E [\|\nabla f_\xi(\vx_j) - \nabla f_\xi(\vx_j)\|^2|\vx_j]\leq \sigma_g^2$ and that $\xi_j$ is drawn independently from the past, then, by Jensen's inequality, we have:
$\E\|\vn_0\|^{3/2}\leq \sigma_{g,0}^{3/2}$ and 
\begin{align*}
    \E\|\sum_{i=0}^{t-1} (1-\alpha)^i\vn_{t-i}\|^{3/2}&\leq  \Big(\E\|\sum_{i=0}^{t-1} (1-\alpha)^i\vn_{t-i}\|^2\Big)^{3/4}\\
    &=\Big(\sum_{i=0}^{t-1} (1-\alpha)^{2i}\E\|\vn_{t-i}\|^2\Big)^{3/4}\\
    &\leq \Big(\sum_{i=0}^{t-1} (1-\alpha)^{2i}\sigma_g^2\Big)^{3/4}\\
    &\leq \Big(\sum_{i=0}^{t-1} (1-\alpha)^{i}\Big)^{3/4}\sigma_g^{3/2}\\
    &\leq \Big(\frac{1}{\alpha}\Big)^{3/4}\sigma_g^{3/2}.\\
\end{align*}

Which means:
$$\E\|\veps_t\|^{3/2}\leq \sqrt{3}\Big((1 - \alpha)^{3t/2}\sigma_{g,0}^{3/2} + \frac{(1-\alpha)^{3/2} L^{3/2}}{\alpha^{2}}\sum_{i=0}^{t-1} (1-\alpha)^i \E r_{t-i}^{3} + \alpha^{3/4}\sigma_g^{3/2}\Big)$$

Finally, 
\begin{align*}
    \frac{1}{T}\sum_{t=0}^{T-1}\E\|\veps_t\|^{3/2}&\leq \frac{\E\|\veps_0\|^{3/2}}{T} +\frac{\sqrt{3}}{T} \sum_{t=1}^{T-1}\Big((1 - \alpha)^{3t/2}\sigma_{g,0}^{3/2} + \frac{(1-\alpha)^{3/2} L^{3/2}}{\alpha^{2}}\sum_{i=0}^{t-1} (1-\alpha)^i \E r_{t-i}^{3} + \alpha^{3/4}\sigma_g^{3/2}\Big)\\
    &\leq \frac{1}{ T} \sigma_{g,0}^{3/2} + \frac{\sqrt{3}}{\alpha T} \sigma_{g,0}^{3/2} + \sqrt{3} \alpha^{3/4}\sigma_g^{3/2} + \frac{\sqrt{3}}{T}\frac{(1-\alpha)^{3/2} L^{3/2}}{\alpha^{2}}\sum_{t=0}^{T-1} \sum_{i=0}^{t-1} (1-\alpha)^i \E r_{t-i}^{3} \\
    &\leq \frac{3}{\alpha T} \sigma_{g,0}^{3/2} + \sqrt{3} \alpha^{3/4}\sigma_g^{3/2} + \frac{\sqrt{3}(1-\alpha)^{3/2} L^{3/2}}{\alpha^{3}}\frac{1}{T}\sum_{t=0}^{T-1}\E r_{t}^{3}
\end{align*}
Which finishes the proof.
\end{proof}
\textbf{Bounding the Hessian momentum error.}

We will do the same for the sequence $(\mH_t)$. Let's define $\mSigma_t = \mH_t - \nabla^2 f(\vx_t)$ and $\mZ(\va,\vb) = \nabla^2 f(\va) - \nabla^2 f(\vb)$ and note that by \eqref{AppLipHess}, we have $\|\mZ(\va,\vb)\|\leq L\|\va - \vb\|$.

\begin{framedlemma}
We show :
\begin{align*}
    \frac{1}{T}\sum_{t=0}^{T-1}\E\|\mSigma_t\|^3\leq \frac{10}{\beta T}\sigma_{h,0}^3 + \beta^{3/2} \Tilde{\sigma}_h^3(\beta) + \frac{9(1 - \beta)L^3}{\beta^3}\frac{1}{T}\sum_{t=0}^{T-1} \E r_{t}^3
\end{align*}

where $\Tilde{\sigma}_h^3(\beta)= 8(er)^{3/2} \sigma_h^3 + 16 (er)^3\beta^{3/2} \delta_h^3$ and $r= 2 \log(d)$.
\end{framedlemma}

\begin{proof}
We have \begin{align*}
    \mSigma_t &= (1 - \beta) \mH_{t-1} + \beta \nabla^2 f_{\xi_t}(\vx_t) - \nabla^2 f(\vx_t)\\
    &= (1 - \beta) \Big(\mH_{t-1}-\nabla^2 f(\vx_t)\Big) + \beta \Big(\nabla^2 f_{\xi_t}(\vx_t) - \nabla^2 f(\vx_t)\Big)
    \\
    &= (1 - \beta) \Big(\mSigma_{t-1}+ \nabla^2 f(\vx_{t-1}) -\nabla^2 f(\vx_t)\Big) + \beta \Big(\nabla^2 f_{\xi_t}(\vx_t) - \nabla^2 f(\vx_t)\Big)
     \\
    &= (1 - \beta)\mSigma_{t-1}+ (1 - \beta)\mZ(\vx_{t-1},\vx_t) + \beta \mN_t
\end{align*}
Where we defined the $\mN_t := \nabla^2 f_{\xi_t}(\vx_t) - \nabla^2 f(\vx_t)$.

Thus  $$\mSigma_t = (1 - \beta)^t \mSigma_{0}+ (1 - \beta)\sum_{i=0}^{t-1}\mZ(\vx_{t-i-1},\vx_{t-i}) + \beta \sum_{i=0}^{t-1}\mN_{t-i}$$

we take the norm now:
\begin{align*}
    \|\mSigma_t\| &\leq (1 - \beta)^t \|\mSigma_{0}\|+ (1 - \beta)\sum_{i=0}^{t-1}(1 - \beta)^i\|\mZ(\vx_{t-i-1},\vx_{t-i})\| + \beta \|\sum_{i=0}^{t-1}(1 - \beta)^i\mN_{t-i}\|\\
    &\leq (1 - \beta)^t \|\mSigma_{0}\|+ (1 - \beta)L\sum_{i=0}^{t-1}(1 - \beta)^ir_{t-i} + \beta \|\sum_{i=0}^{t-1}(1 - \beta)^i\mN_{t-i}\|
\end{align*}

We use the convexity of $x\mapsto x^3$ and get :
\begin{align*}
    \|\mSigma_t\|^3 &\leq 9\Big((1 - \beta)^{3t} \|\mSigma_{0}\|^3+ (1 - \beta)L^3(\sum_{i=0}^{t-1}(1 - \beta)^ir_{t-i})^3 + \beta^3 \|\sum_{i=0}^{t-1}(1 - \beta)^i\mN_{t-i}\|^3\Big)\\
    &\leq 9\Big((1 - \beta)^{3t} \|\mN_{0}\|^3+ \frac{(1 - \beta)L^3}{\beta^2}\sum_{i=0}^{t-1}(1 - \beta)^i r_{t-i}^3 + \beta^3 \|\sum_{i=0}^{t-1}(1 - \beta)^i\mN_{t-i}\|^3\Big)
\end{align*}

We take the expectation and get 
\begin{equation}\label{HN0}
    \E\|\mSigma_t\|^3 \leq 9\Big((1 - \beta)^{3t} \E\|\mN_{0}\|^3+ \frac{(1 - \beta)L^3}{\beta^2}\sum_{i=0}^{t-1}(1 - \beta)^i \E r_{t-i}^3 + \beta^3 \E\|\sum_{i=0}^{t-1}(1 - \beta)^i\mN_{t-i}\|^3\Big)
\end{equation}

We make the following assumptions about the Hessian noises $\mN_t$:
\begin{itemize}
    \item the noise is centered: $\E\mN_t = 0$;
    \item  the noise has a bounded variance: there exists $\sigma_h\geq 0$ such that $\E\|\mN_t\|^2\leq \sigma_h$;
    \item the noise is bounded almost surely, i.e. there exists $\delta_h\geq \sigma_h$ such that $\|\mN_t\|\leq \delta_h$ with probability one.
\end{itemize}
We also assume that the noise at the start $\E\|\mN_0\|^2\leq \sigma_{h,0}^2$ where $\sigma_{h,0}$ is potentially much smaller than $\sigma_{h}$ (e.g., by using a large batch at the start).

Note that all these noises are independent from each other conditioned on the past.

Under these assumptions we can show that:
$\E\|\mN_{0}\|^3\leq \sigma_{h,0}^3$ and using the following Lemma ~\ref{lyapMatr},
we have for $q = 3$ and $r = 2 log(d)$ where $d$ is the dimension of the problem:

\begin{align*}
    \Big(\E\|\sum_{i=0}^{t-1}(1 - \beta)^i\mN_{t-i}\|^3\big)^{1/3}&\leq 2\sqrt{er} \|\big(\sum_{i=0}^{t-1}(1 - \beta)^i \E[\mN_{t-i}^2]\big)^{1/2}\| + 4 er \E[\max_i \|(1 - \beta)^i\mN_{t-i}\|^q]^{1/q}\\
    &\leq 2\sqrt{er} \|\big(\sum_{i=0}^{t-1}(1 - \beta)^i \E[\mN_{t-i}^2]\big)\|^{1/2} + 4 er \E[\max_i \|(1 - \beta)^i\mN_{t-i}\|^q]^{1/q}\\
    &\leq 2\sqrt{er} \big(\sum_{i=0}^{t-1}(1 - \beta)^i \big)^{1/2} \sigma_h + 4 er \delta_h\\
    &\leq 2\sqrt{\frac{er}{\beta}} \sigma_h + 4 er \delta_h.
\end{align*}
We conclude that :
\begin{align*}
    \beta^3\E\|\sum_{i=0}^{t-1}(1 - \beta)^i\mN_{t-i}\|^3&\beta^3\leq \Big( 2\sqrt{\frac{er}{\beta}} \sigma_h + 4 er \delta_h\Big)^3\\
    &\leq 8(er\beta)^{3/2} \sigma_h^3 + 16 (er\beta)^3 \delta_h^3\\
    &= \beta^{3/2} \Big(8(er)^{3/2} \sigma_h^3 + 16 (er)^3\beta^{3/2} \delta_h^3\Big)\\
    &:= \beta^{3/2} \Tilde{\sigma}_h^3(\beta).
\end{align*}

Replacing this in \eqref{HN0}, we get:
$$\E\|\mSigma_t\|^3 \leq 9\Big((1 - \beta)^{3t}\sigma_{h,0}^3 + \beta^{3/2} \Tilde{\sigma}_h^3(\beta) + \frac{(1 - \beta)L^3}{\beta^2}\sum_{i=0}^{t-1}(1 - \beta)^i \E r_{t-i}^3\Big).$$

Finally:
\begin{align*}
    \frac{1}{T}\sum_{t=0}^{T-1}\E\|\mSigma_t\|^3&\leq \frac{\E\|\mSigma_0\|^3}{T} + \frac{9}{T}\sum_{t=1}^{T-1}\Big((1 - \beta)^{3t}\sigma_{h,0}^3 + \beta^{3/2} \Tilde{\sigma}_h^3(\beta) + \frac{(1 - \beta)L^3}{\beta^2}\sum_{i=0}^{t-1}(1 - \beta)^i \E r_{t-i}^3\Big)\\
    &\leq \frac{1}{ T}\sigma_{h,0}^3 + \frac{9}{\beta T}\sigma_{h,0}^3 + \beta^{3/2} \Tilde{\sigma}_h^3(\beta) + \frac{9}{T}\frac{(1 - \beta)L^3}{\beta^2}\sum_{t=0}^{T-1} \sum_{i=0}^{t-1}(1 - \beta)^i \E r_{t-i}^3\\
    &\leq \frac{10}{\beta T}\sigma_{h,0}^3 + \beta^{3/2} \Tilde{\sigma}_h^3(\beta) + \frac{9(1 - \beta)L^3}{\beta^3}\frac{1}{T}\sum_{t=0}^{T-1} \E r_{t}^3.
\end{align*}
\end{proof}
\subsection{Nonconvex convergence}

Using Theorem~\ref{ThOneStep}, we have:
\begin{align*}
    f(\vx_{t}) - f(\vx_{t+1}) &\geq
 \frac{1}{1008\sqrt{M}}\mu_M(\vx_{t+1}) 
+ \frac{M \| \vx_{t} - \vx_{t+1}\|^3 }{72} - \frac{4 \|\nabla f(\vx_{t}) - \vg_t\|^{3/2}}{\sqrt{M}} - \frac{73 \| \nabla^2 f(\vx_t) - \mH_t \|^3}{M^2}\\
&\geq
 \frac{1}{1008\sqrt{M}}\mu_M(\vx_{t+1}) 
+ \frac{M r_{t}^3 }{72} - \frac{4 \|\veps_t\|^{3/2}}{\sqrt{M}} - \frac{73 \| \mSigma_t \|^3}{M^2}.
\end{align*}

Reorganizing the inequality and taking the expectation we get:

$$\frac{1}{1008}\E\mu_M(\vx_{t+1})\leq \sqrt{M}\E(f(\vx_{t}) - f(\vx_{t+1})) - \frac{M^{3/2} \E r_{t}^3 }{72} + 4 \E\|\veps_t\|^{3/2} + \frac{73 \E\| \mSigma_t \|^3}{M^{3/2}}.$$
Summing over $t=0, T-1$ and dividing by $T$, we get:
\begin{align*}
    \frac{1}{1008} \frac{1}{T}\sum_{t=0}^{T-1}\E\mu_M(\vx_{t+1})\leq \frac{\sqrt{M}F_0}{T} - \frac{M^{3/2}}{72} \E \frac{1}{T}\sum_{t=0}^{T-1}r_{t}^3 + 4 \frac{1}{T}\sum_{t=0}^{T-1}\E\|\veps_t\|^{3/2} + \frac{73 }{M^{3/2}}\frac{1}{T}\sum_{t=0}^{T-1}\E\| \mSigma_t \|^3,
\end{align*}

where $F_0 := f(\vx_0) - \min_\vx f(\vx)$ and we assume $F_0<+\infty$.

We have already shown that :
\begin{align*}
    \frac{1}{T}\sum_{t=0}^{T-1}\E\|\veps_t\|^{3/2}
    &\leq \frac{3}{\alpha T} \sigma_{g,0}^{3/2} + \sqrt{3} \alpha^{3/4}\sigma_g^{3/2} + \frac{\sqrt{3}(1-\alpha)^{3/2} L^{3/2}}{\alpha^{3}}\frac{1}{T}\sum_{t=0}^{T-1}\E r_{t}^{3}
\end{align*}

and 
\begin{align*}
    \frac{1}{T}\sum_{t=0}^{T-1}\E\|\mSigma_t\|^3
    &\leq \frac{10}{\beta T}\sigma_{h,0}^3 + \beta^{3/2} \Tilde{\sigma}_h^3(\beta) + \frac{9(1 - \beta)L^3}{\beta^3}\frac{1}{T}\sum_{t=0}^{T-1} \E r_{t}^3.
\end{align*}

Thus we get:

\begin{align*}
    \frac{1}{1008} \frac{1}{T}\sum_{t=0}^{T-1}\E\mu_M(\vx_{t+1})&\leq \frac{\sqrt{M}F_0}{T} - \frac{M^{3/2}}{72} \E \frac{1}{T}\sum_{t=0}^{T-1}r_{t}^3 + 4 \frac{1}{T}\sum_{t=0}^{T-1}\E\|\veps_t\|^{3/2} + \frac{73 }{M^{3/2}}\frac{1}{T}\sum_{t=0}^{T-1}\E\| \mSigma_t \|^3\\
    &\leq \frac{\sqrt{M}F_0}{T} +  \frac{12}{\alpha T} \sigma_{g,0}^{3/2} + 4\sqrt{3} \alpha^{3/4}\sigma_g^{3/2} + \frac{730 }{M^{3/2}\beta T}\sigma_{h,0}^3 + \frac{73 }{M^{3/2}}\beta^{3/2} \Tilde{\sigma}_h^3(\beta)\\& + \Big(\frac{4\sqrt{3} L^{3/2}}{\alpha^{3}} + \frac{73 }{M^{3/2}}\frac{9L^3}{\beta^3} - \frac{M^{3/2}}{72}\Big)\E\frac{1}{T}\sum_{t=0}^{T-1}r_{t}^3.
\end{align*}

We choose $\alpha,\beta \in(0,1]$ such that:
    \begin{equation}\label{Condition_app}
        \frac{4\sqrt{3} L^{3/2}}{\alpha^{3}} + \frac{73 }{M^{3/2}}\frac{9L^3}{\beta^3} - \frac{M^{3/2}}{72}\leq 0\;,
    \end{equation}
and to deal with initial gradient and Hessian noises, we also take:

$\alpha \in \argmin \frac{12}{\alpha T} \sigma_{g,0}^{3/2} + 4\sqrt{3} \alpha^{3/4}\sigma_g^{3/2}\approx \frac{a_g^{6/7}}{T^{4/7}}$ and $\beta \in \argmin \frac{730 }{M^{3/2}\beta T}\sigma_{h,0}^3 + \frac{73 }{M^{3/2}}\beta^{3/2} \sigma_h^3\approx  \frac{a_h^{6/5}}{T^{2/5}}$.

Where, for simplicity, we denoted $a_g = \frac{\sigma_{g,0}}{\sigma_{g}}\leq 1$ and $a_h = \frac{\sigma_{h,0}}{\sigma_{h}}\leq 1$. 

All in all, we choose 
\begin{equation}\label{ChoiceOfParams:App}
\ba{rcl}
    \alpha = \max\Big(\frac{a_g^{6/7}}{T^{4/7}},10\sqrt{\frac{L}{M}}\Big)\, , \,
    \beta = \max\Big(\frac{a_h^{6/5}}{T^{2/5}},46\frac{L}{M}\Big)\,.
\ea
\end{equation}
To ensure $\alpha, \beta \leq 1$ it suffices to take $M\geq 100 L$.

Finally, we get:

\begin{align*}
    \frac{1}{1008} \frac{1}{T}\sum_{t=0}^{T-1}\E\mu_M(\vx_{t+1})
    &\leq \frac{\sqrt{M}F_0}{T}  + 40\frac{L^{3/8}}{M^{3/8}}\sigma_g^{3/2} + \frac{73 (46)^{3/2}L^{3/2} }{M^{3}} \Tilde{\sigma}_h^3(\beta) + 20 \frac{a_{g}^{9/14}\sigma_{g}^{3/2}}{T^{3/7}} + 803\frac{a_h^{9/5}\sigma_{h}^{3}}{M^{3/2} T^{3/5}},
\end{align*}

where $\Tilde{\sigma}_h^3(\beta)= 8(er)^{3/2} \sigma_h^3 + 16 (46)^{3/2}(er)^3\frac{L^{3/2}}{M^{3/2}} \delta_h^3$ (remember that $r= 2 \log(d)$).

Compare this to the bound we get without using momentum:
$$\cO\Big(\frac{\sqrt{M}F_0}{T}  + \sigma_g^{3/2} + \frac{1 }{M^{3/2}} \sigma_h^3\Big).$$

We clearly depend better on the regularization parameter $M$ in both the gradient and Hessian noise terms.

Going back to the bound 
\begin{align*}
    \frac{1}{1008} \frac{1}{T}\sum_{t=0}^{T-1}\E\mu_M(\vx_{t+1})
    &\leq \frac{\sqrt{M}F_0}{T}  + 40\frac{L^{3/8}}{M^{3/8}}\sigma_g^{3/2} + \frac{73 (46)^{3/2}L^{3/2} }{M^{3}} \Tilde{\sigma}_h^3(\beta) + 20 \frac{a_{g}^{9/14}\sigma_{g}^{3/2}}{T^{3/7}} + 803\frac{a_h^{9/5}\sigma_{h}^{3}}{M^{3/2} T^{3/5}}.
\end{align*}

By choosing $$M = \max\Big(100 L, \frac{L^{3/7}\sigma_g^{12/7}}{F_0^{8/7}}T^{8/7}, \frac{L^{3/7}\Tilde{\sigma}_h(T)^{3/7}}{F_0^{2/7}}T^{2/7}, \frac{a_h^{9/10} \sigma_h^{3/2} T^{1/5}}{\sqrt{F_0}}\Big)$$
for $\Tilde{\sigma}_h^3(T)=8(er)^{3/2} \sigma_h^3 + 16 (46)^{3/2}(er)^3\frac{L^{6/7}\Tilde{F_0}^{3/7}}{\sigma_h^{9/7}T^{3/7}} \delta_h^3$, we get:

\begin{align*}
    \frac{1}{1008} \frac{1}{T}\sum_{t=0}^{T-1}\E\mu_M(\vx_{t+1})
    &\leq \frac{\sqrt{L}F_0}{T}  + \frac{40L^{3/14}F_0^{3/7}\sigma_g^{6/7} + 20 a_{g}^{9/14}\sigma_{g}^{3/2}}{T^{3/7}} + 73 (46)^{3/2} \frac{L^{3/14}F_0^{6/7}\Tilde{\sigma}_h(T)^{3/7}}{T^{6/7}} + 804\frac{a_h^{9/20}\sigma_{h}^{3/4}F_0^{1/2}}{ T^{9/10}}.
\end{align*}

\subsection{Convex case}

We consider here the slightly more general stochastic composite optimization problem 
\begin{equation}
    \min_\vx f(\vx):=\underbrace{\E_{\xi\sim \mathcal{P}}\Big[f_\xi(\vx)\Big]}_{:=\Tilde{f}(\vx)} + \psi(\vx)\;,
\end{equation}

where $\psi$ is a simple, convex, potentially non-smooth function and $f$ is convex. Simple here means that we can efficiently solve the following general subproblem:
\begin{equation}\label{xx+cv}
\ba{rcl}
\!\!\!\!\!\!
    \vx_{t+1} & \in & \argmin_{\vy \in \R^d}
    \Bigl\{\, \Omega_{M,\vg_t,\mH_t}(\vy,\vx_{t})  + \psi(\vy)\Bigr\}\,.
\ea
\end{equation}
For example, $\psi$ can be the indicator function of a bounded convex set, in this case, solving \eqref{xx+cv} has the same cost as solving the cubic subproblem \eqref{xx+}.

The optimality condition of \eqref{xx+cv} guarantees that there exists a subgradient of $\psi$ that we denote as $\psi^\prime(\vx_{t+1}) \in \partial \psi(\vx_{t + 1})$ such that $$\vg_t + \mH_t (\vx_{t+1} - \vx_t) + \frac{M}{2}r_t^2 + \psi^\prime(\vx_{t+1}) = \vzero.$$
 We define $\nabla f(\vx_{t+1}) = \nabla \Tilde{f}(\vx_{t+1}) + \psi^\prime(\vx_{t+1})$.

We assume that the diameter of the domain of $f$ is bounded by $D$.
This means that :
$$\forall\vx\in dom(f)\quad:\quad f(\vx) - f^\star \leq D \|\nabla f(\vx)\|\;.$$

Following the same lines as  Theorem C.1 from \citep{chayti2024unified}, it is easy to show the same inequality as the one we had before for the progress: 
\begin{align*}
    f(\vx_{t}) - f(\vx_{t+1})
&\geq
 \frac{1}{1008\sqrt{M}}\|\nabla f(\vx_{t+1})\|^{3/2} 
+ \frac{M r_{t}^3 }{72} - \frac{4 \|\veps_t\|^{3/2}}{\sqrt{M}} - \frac{73 \| \mSigma_t \|^3}{M^2}\\&\geq
 \frac{1}{1008\sqrt{M}D^{3/2}}(f(\vx_{t+1} - f^\star)^{3/2} 
+ \frac{M r_{t}^3 }{72} - \frac{4 \|\veps_t\|^{3/2}}{\sqrt{M}} - \frac{73 \| \mSigma_t \|^3}{M^2}.
\end{align*}
Let us define $F_t = \E f(\vx_{t}) - f^\star $, using Jensen's inequality we get 
\begin{equation}\label{convex}
    F_t - F_{t+1}
\geq
 \frac{1}{1008\sqrt{M}D^{3/2}} F_{t+1}^{3/2} 
+ \underbrace{\E\Big[\frac{M r_{t}^3 }{72} - \frac{4 \|\veps_t\|^{3/2}}{\sqrt{M}} - \frac{73 \| \mSigma_t \|^3}{M^2}\Big]}_{:= - a_t}.
\end{equation}

From the non-convex case, we have that 

$$ \frac{\sqrt{M}}{T} \sum_{t=0}^{T-1}a_t \leq \frac{12}{\alpha T} \sigma_{g,0}^{3/2} + 4\sqrt{3} \alpha^{3/4}\sigma_g^{3/2} + \frac{730 }{M^{3/2}\beta T}\sigma_{h,0}^3 + \frac{73 }{M^{3/2}}\beta^{3/2} \Tilde{\sigma}_h^3(\beta)$$
under the condition that 
$$\frac{4\sqrt{3} L^{3/2}}{\alpha^{3}} + \frac{73 }{M^{3/2}}\frac{9L^3}{\beta^3} - \frac{M^{3/2}}{72}\leq 0\;.$$
Thus using the same choices of $\alpha,\beta$ as in the non-convex case, we get:
\begin{equation}\label{a_t bound}
     \frac{\sqrt{M}}{T} \sum_{t=0}^{T-1}a_t \leq 20 \frac{a_{g}^{9/14}\sigma_{g}^{3/2}}{T^{3/7}} + 803\frac{a_h^{9/5}\sigma_{h}^{3}}{M^{3/2} T^{3/5}} + 40\frac{L^{3/8}}{M^{3/8}}\sigma_g^{3/2} + \frac{73 (46)^{3/2}L^{3/2} }{M^{3}} \Tilde{\sigma}_h^3(\beta)\;,
\end{equation}

Going back to \eqref{convex}, we have:
$$F_t - F_{t+1}
\geq
 C F_{t+1}^{3/2} 
- a_t\;,$$
where we denoted $C:=\frac{1}{1008\sqrt{M}D^{3/2}}$.

Consequently:
\begin{align*}
    \frac{1}{\sqrt{F_{t+1}}} - \frac{1}{\sqrt{F_{t}}} &=\frac{\sqrt{F_t} - \sqrt{F_{t+1}}}{\sqrt{F_t F_{t+1}}}\\
    &\geq\frac{F_t - F_{t+1}}{2 F_t\sqrt{ F_{t+1}}}.
    \intertext{where we used the concavity of $x\mapsto \sqrt{x}$ which gives $\sqrt{a} - \sqrt{b}\leq \frac{1}{2\sqrt{b}}(a - b)$, and took $a= F_{t+1}, b=F_t$.}
    \intertext{Thus:}
    \frac{1}{\sqrt{F_{t+1}}} - \frac{1}{\sqrt{F_{t}}} 
    &\geq\frac{C F_{t+1}^{3/2} 
- a_t}{2F_{t}\sqrt{F_{t+1}}} = \frac{C}{2}\frac{F_{t+1}}{F_t} - \frac{a_t}{2F_{t}\sqrt{F_{t+1}}}.
\intertext{Let $\varepsilon = \underset{0\leq t\leq T-1}{\min} F_t $. Thus, for all $t\leq T - 1$ we have $F_t\geq \varepsilon$, then:}
\frac{1}{\sqrt{\varepsilon}}&\geq \frac{1}{F_T} - \frac{1}{F_0}\\
&\geq \frac{C T}{2} \Big(\frac{1}{T}\sum_{t=0}^{T-1}\frac{F_{t+1}}{F_t}\Big) - \frac{T}{2\varepsilon^{3/2}\sqrt{M}}\Big(\frac{\sqrt{M}}{T}\sum_{t=0}^{T-1}a_t\Big).
\intertext{We note that using the fact that the arithmetic mean is greater than the geometric mean, we have:}
\frac{1}{T}\sum_{t=0}^{T-1}\frac{F_{t+1}}{F_t}&\geq \Big(\prod_{t=0}^{T-1}\frac{F_{t+1}}{F_t}\Big)^{1/T} = \Big(\frac{F_{T}}{F_0}\Big)^{1/T}\geq \Big(\frac{\varepsilon}{F_0}\Big)^{1/T}\geq 1+\frac{1}{T} \log(\frac{\varepsilon}{F_0}),
\intertext{which leads to}
\frac{1}{\sqrt{\varepsilon}}&\geq \frac{C T}{2} \Big(1 - \frac{1}{T} \log(\frac{F_0}{\varepsilon})\Big) - \frac{T}{2\varepsilon^{3/2}\sqrt{M}}\Big(\frac{\sqrt{M}}{T}\sum_{t=0}^{T-1}a_t\Big).
\intertext{Remember \eqref{a_t bound}, then:}
\frac{1}{\sqrt{\varepsilon}}&\geq \frac{C T}{2} \Big(1 - \frac{1}{T} \log(\frac{F_0}{\varepsilon})\Big) - \frac{T}{2\varepsilon^{3/2}\sqrt{M}}\mathcal{A}(M),
\intertext{where we denoted $\mathcal{A}(M) := 20 \frac{a_{g}^{9/14}\sigma_{g}^{3/2}}{T^{3/7}} + 803\frac{a_h^{9/5}\sigma_{h}^{3}}{M^{3/2} T^{3/5}} +40\frac{L^{3/8}}{M^{3/8}}\sigma_g^{3/2} + \frac{73 (46)^{3/2}L^{3/2} }{M^{3}} \Tilde{\sigma}_h^3(\beta)\;.$ We conclude that:}
\frac{1}{\sqrt{\varepsilon}}& + \frac{C}{2}\log(\frac{F_0}{\varepsilon}) + \frac{T}{2\varepsilon^{3/2}\sqrt{M}}\mathcal{A}(M)\geq \frac{C}{2}T,
\intertext{We invert the previous inequality and get: }
\varepsilon& \leq \frac{36}{C^2 T^2} + F_0 e^{-T/3} + \Big(\frac{\mathcal{A}(M)}{C\sqrt{M}}\Big)^{2/3}
\intertext{We conclude that:}
\underset{0\leq t\leq T-1}{\min} F_t &\leq \cO\Big( \frac{M D^3}{T^2} + \frac{L D \Tilde{\sigma}_h^2}{M^2} + \frac{L^{1/4}}{M^{1/4}}D\sigma_g + F_0 e^{-T/3} + \frac{a_g^{3/7} D \sigma_g}{T^{2/7}} + \frac{a_h^{6/5} D \sigma_h^2}{MT^{2/5}} \Big)
\intertext{We choose $M = \max\Big(100 L,\frac{L^{1/5} \sigma_g^{4/5} T^{8/5}}{D^{8/5}},\frac{L^{1/3}\Tilde{\sigma}_h^{2/3}}{D^{2/3}}, \frac{a_h^{3/5} \sigma_h T^{4/5}}{D}\Big)$.}
\end{align*}
By replacing $M$ with its value we arrive at the following:
$$\min_{0\leq t\leq T-1} F_t=\cO\Big( \frac{L D^3}{T^2} + \frac{L^{1/3}D^{7/3}\Tilde{\sigma}_h^{2/3}}{T^{4/3}} + \frac{L^{1/5}D^{7/5}\sigma_g^{4/5}}{T^{2/5}} + F_0 e^{-T/3} + \frac{a_h^{3/5} D^{2}\sigma_h}{T^{6/5}} + \frac{a_g^{3/7} D \sigma_g}{T^{2/7}}  \Big).$$

%We notice that this upper bound on $T$ is badly affected by the terms proportional to $\sigma_{g,0}^{3/2}$ and $\sigma_{h,0}^{3}$, we remind that $\sigma_{g,0}$ ($\sigma_{h,0}$) represents the variance of the first gradient (Hessian resp.) sample. One simple fix that can improve this rate is to use a large batch to have a good first estimate of the gradient and the Hessian; let this batch have size $b$, then we will have:
%$$\sigma_{g,0} = \sigma_g/\sqrt{b}\; ,$$
%we can, for example, choose $b = \frac{T}{10}$ or any constant other than $10$, then this will only change the rate by a multiplicative factor of $\frac{11}{10}$.
%By doing this we can show that:
%$$\frac{\sqrt{MD^3}\sigma_{g,0}^{3/2}}{\varepsilon^{3/2}} = \Tilde{\cO}\Big(\frac{L^{1/8}D^{3/8}}{\varepsilon^{9/8}} + \frac{L^{1/16}D^{7/16}\Tilde{\sigma}_h}{\varepsilon^{19/16}} + \frac{L^{1/8}D^{7/8}\sigma_g^{1/2}}{\varepsilon^{13/8}} \Big)\sigma_g^{3/2}$$
%and $$\frac{D^{3/2} \sigma_{h,0}^{3}}{L\sqrt{M}\varepsilon^{3/2}} = \Tilde{\cO}\Big( \frac{ \sigma_{h}^{3}}{L D^{3/4}M^2\varepsilon^{3/4}}\Big).$$
%Note that all powers of $\varepsilon$ are smaller than $\frac{1}{2}$.

%Thus, in this case, 

%\begin{multline*}
%$$
%T=\cO\Big(\log(\frac{F_0}{\varepsilon}) + \frac{\sqrt{LD^3}}{\sqrt{\varepsilon}} + \frac{D^{7/4}L^{1/4}\Tilde{\sigma}_h^{1/2}}{\varepsilon^{3/4}} + \frac{D^{5/2}L^{1/4}\sigma_g}{\varepsilon^{5/2}}\Big).
%$$
%\end{multline*}

\section{Second-order momentum for the gradients} \label{sec:SOM}
We will now consider the following momentum:

\begin{equation}\label{SOM}
    \vg_t = (1 - \alpha) \Big(\vg_{t-1} + \mH_t(\vx_t - \vx_{t-1})\Big) + \alpha \nabla f_{\xi_t}(\vx_t ), \;\text{for}\;t\geq 1\;\text{and}\; \vg_0 = \nabla f_{\xi_0}(\vx_0),
\end{equation}

where $\mH_t$ is the Hessian momentum defined as
\begin{equation}
    \mH_t = (1 - \beta) \mH_{t-1} + \beta \nabla^2 f_{\xi_t}(\vx_t), \;\text{for}\;t\geq 1\;\text{and}\; \vg_0 = \nabla^2 f_{\xi_0}(\vx_0)\;,
\end{equation}
for $\alpha,\beta\in(0,1]$.

\begin{framedlemma} For the update in \eqref{SOM}, we have for any $\gamma>0$:
  \begin{multline*}
      \frac{1}{2T}\sum_{t=0}^{T-1}\E\|\veps_t\|^{3/2}\leq\frac{\sigma_{g,0}^{3/2}}{\alpha T} + \frac{(1-\beta)\sigma_{h,0}^3}{\gamma^{3/2}\alpha^{3/2} \beta T} +  \alpha^{3/4}\sigma_g^{3/2} + \frac{\beta^{3/2}}{\gamma^{3/2}\alpha^{3/2}} \Tilde{\sigma}_h^3(\beta) \\+(1 - \alpha) \Big(\frac{L^{3/2}}{\alpha^{3/2}} + \frac{\gamma^{3/2}}{\alpha^{3/2}} + \frac{1 - \beta}{2^{3/2}\gamma^{3/2}\alpha^{3/2}\beta^3}\Big)\frac{1}{T}\sum_{t=0}^{T-1}\E r^3_{t}.
  \end{multline*}
\end{framedlemma}
\begin{proof}
    
As before, let $\veps_t = \vg_t - \nabla f(\vx_t)$, then we have:
\begin{align*}
    \veps_t &= (1 - \alpha) \Big(\vg_{t-1} + \mH_t(\vx_t - \vx_{t-1})\Big) + \alpha \nabla f_{\xi_t}(\vx_t ) - \nabla f(\vx_t)\\
    &= (1 - \alpha) \Big(\veps_{t-1} + \nabla f(\vx_{t-1}) - \nabla f(\vx_t) + \mH_t(\vx_t - \vx_{t-1})\Big) + \alpha (\nabla f_{\xi_t}(\vx_t ) - \nabla f(\vx_t))\\
    &= (1 - \alpha) \Big(\veps_{t-1} + \vz(\vx_{t-1},\vx_{t}) + [\mH_t - \nabla^2f(\vx_t)](\vx_t - \vx_{t-1})\Big) + \alpha (\nabla f_{\xi_t}(\vx_t ) - \nabla f(\vx_t))\\
    &= (1 - \alpha) \Big(\veps_{t-1} + \vz(\vx_{t-1},\vx_{t}) + \mSigma_t(\vx_t - \vx_{t-1})\Big) + \alpha \vn_t.
    \intertext{Thus}
    \veps_t &= (1 - \alpha)^t\veps_0 + (1 - \alpha)\sum_{i=0}^{t-1}(1 - \alpha)^i\Big(\vz(\vx_{t-i-1},\vx_{t-i}) + \mSigma_{t-i}(\vx_{t-i} - \vx_{t-i-1})\Big) + \alpha\sum_{i=0}^{t-1}(1 - \alpha)^i\vn_{t-i}.
    \intertext{By taking the norm, we get:}
    \|\veps_t\|&\leq(1 - \alpha)^t\|\vn_0\| + (1 - \alpha)\sum_{i=0}^{t-1}(1 - \alpha)^i\Big(\frac{L}{2}r^2_{t-i} + \|\mSigma_{t-i}\|r_{t-i}\Big) + \alpha\|\sum_{i=0}^{t-1}(1 - \alpha)^i\vn_{t-i}\|\\
    &\leq(1 - \alpha)^t\|\vn_0\| + (1 - \alpha)\frac{L+\gamma}{2}\sum_{i=0}^{t-1}(1 - \alpha)^ir^2_{t-i} + (1 - \alpha)\frac{1}{2\gamma}\sum_{i=0}^{t-1}(1 - \alpha)^i\|\mSigma_{t-i}\|^2 + \alpha\|\sum_{i=0}^{t-1}(1 - \alpha)^i\vn_{t-i}\|,
    \intertext{where $\gamma > 0$ is to be chosen latter. Now we take the power $3/2$ and get:}
    \frac{1}{2}\|\veps_t\|^{3/2}&\leq(1 - \alpha)^t\|\vn_0\|^{3/2} + (1 - \alpha)\Big(\frac{L+\gamma}{2}\sum_{i=0}^{t-1}(1 - \alpha)^ir^2_{t-i}\Big)^{3/2} + (1 - \alpha)\Big(\frac{1}{2\gamma}\sum_{i=0}^{t-1}(1 - \alpha)^i\|\mSigma_{t-i}\|^2\Big)^{3/2} \\&+ \alpha^{3/2}\|\sum_{i=0}^{t-1}(1 - \alpha)^i\vn_{t-i}\|^{3/2}\\
    &\leq(1 - \alpha)^t\|\vn_0\|^{3/2} + (1 - \alpha)\frac{(L+\gamma)^{3/2}}{2^{3/2}\alpha^{1/2}}\sum_{i=0}^{t-1}(1 - \alpha)^ir^3_{t-i} + (1 - \alpha)\frac{1}{2^{3/2}\gamma^{3/2}\alpha^{1/2}}\sum_{i=0}^{t-1}(1 - \alpha)^i\|\mSigma_{t-i}\|^{3} \\&+ \alpha^{3/2}\|\sum_{i=0}^{t-1}(1 - \alpha)^i\vn_{t-i}\|^{3/2}.
    \intertext{By taking the expectation, and using the same reasoning as before, we get:}
    \frac{1}{2}\E\|\veps_t\|^{3/2}&\leq(1 - \alpha)^t\sigma_{g,0}^{3/2} + (1 - \alpha)\frac{(L+\gamma)^{3/2}}{2^{3/2}\alpha^{1/2}}\sum_{i=0}^{t-1}(1 - \alpha)^i\E r^3_{t-i} + \alpha^{3/4}\sigma_g^{3/2}+ (1 - \alpha)\frac{1}{2^{3/2}\gamma^{3/2}\alpha^{1/2}}\sum_{i=0}^{t-1}(1 - \alpha)^i\E\|\mSigma_{t-i}\|^{3}.
\end{align*}
By substituting our bound of $\E\|\mSigma_{t-i}\|^{3} $ from before, we get:
\begin{align*}
    \frac{1}{9}\sum_{i=0}^{t-1}(1 - \alpha)^i\E\|\mSigma_{t-i}\|^{3}&\leq \sum_{i=0}^{t-1}(1 - \alpha)^{i} \Big((1 - \beta)^{t-i}\sigma_{h,0}^3 + \beta^{3/2} \Tilde{\sigma}_h^3(\beta) + \frac{(1 - \beta)L^3}{\beta^2}\sum_{j=0}^{t-i-1}(1 - \beta)^j \E r_{t-i-j}^3\Big)\\
    &\leq \sum_{i=0}^{t-1}(1 - \alpha)^{i} \Big((1 - \beta)^{t-i}\sigma_{h,0}^3 + \beta^{3/2} \Tilde{\sigma}_h^3(\beta) + \frac{(1 - \beta)L^3}{\beta^2}\sum_{j=0}^{t-i-1}(1 - \beta)^j \E r_{t-i-j}^3\Big)\\
    &\leq (1 - \beta)\frac{(1 -\beta)^t -(1 -\alpha)^t }{\alpha - \beta}\sigma_{h,0}^3 + \frac{\beta^{3/2}}{\alpha} \Tilde{\sigma}_h^3(\beta) + \frac{(1 - \beta)L^3}{\beta^2}\sum_{i=0}^{t-1}(1 - \alpha)^{i}\sum_{j=0}^{t-i-1}(1 - \beta)^j \E r_{t-i-j}^3\\
    &\leq (1 - \beta)\frac{(1 -\min(\alpha,\beta))^t }{\max(\alpha,\beta)}\sigma_{h,0}^3 + \frac{\beta^{3/2}}{\alpha} \Tilde{\sigma}_h^3(\beta) + \frac{(1 - \beta)L^3}{\beta^2}\sum_{i=0}^{t-1}\frac{(1 -\min(\alpha,\beta))^i  }{\max(\alpha,\beta)} \E r_{t-i}^3.
\end{align*}
Thus,
\begin{align*}
    \frac{1}{2T}\sum_{t=0}^{T-1}\E\|\veps_t\|^{3/2}&\leq\frac{\sigma_{g,0}^{3/2}}{\alpha T} + \frac{(1-\beta)\sigma_{h,0}^3}{\gamma^{3/2}\alpha^{3/2} \beta T} +  \alpha^{3/4}\sigma_g^{3/2} + \frac{\beta^{3/2}}{\gamma^{3/2}\alpha^{3/2}} \Tilde{\sigma}_h^3(\beta) \\&+(1 - \alpha) \Big(\frac{L^{3/2}}{\alpha^{3/2}} + \frac{\gamma^{3/2}}{\alpha^{3/2}} + \frac{1 - \beta}{2^{3/2}\gamma^{3/2}\alpha^{3/2}\beta^3}\Big)\frac{1}{T}\sum_{t=0}^{T-1}\E r^3_{t}.
\end{align*}
\end{proof}
This means we need to choose $\alpha, \beta$ such that $$\frac{L^{3/2}}{\alpha^{3/2}} + \frac{\gamma^{3/2}}{\alpha^{3/2}} + \frac{(1 - \beta) L^3}{2^{3/2}\gamma^{3/2}\alpha^{3/2}\beta^3} + (1 - \beta)\frac{73 }{M^{3/2}}\frac{9L^3}{\beta^3} - \frac{M^{3/2}}{72} \leq 0.$$

For example, for $\beta = 1$ (smaller choices of $\beta$ will make the middle term very high and are less likely to lead to any improvement), it suffices to take:
$\alpha = 3175 \max(\frac{L}{M},\frac{\gamma}{M})$.

We choose $\gamma$ such that it minimizes $\frac{\gamma^{3/4}}{M^{3/4}}\sigma_g^{3/2} + \frac{1}{\gamma^{3/2}\alpha^{3/2}} \sigma_h^3$ and $\gamma\leq M/3175$ to ensure $\alpha$ can be made $\leq 1$. We choose $\gamma = \min (M,\frac{M^{3/5}\sigma_h^{4/5}}{\sigma_g^{2/5}})$. Note that when $\gamma = M/3175$, we have $\sigma_g \leq 3175^{5/2}\frac{\sigma_h^2}{M}$.

For this choice of $\gamma$ we get:
$\frac{\gamma^{3/4}}{M^{3/4}}\sigma_g^{3/2} + \frac{1}{\gamma^{3/2}\alpha^{3/2}} \sigma_h^3 \leq \frac{2\sigma_h^3}{M^{3/2}} + \frac{2 \sigma_g^{9/10}\sigma_h^{3/5}}{M^{3/10}}$.

We also choose $M\geq 3175 L$ to ensure $\alpha \leq 1.$

 Thus 
 \begin{align*}
    \frac{1}{1008} \frac{1}{T}\sum_{t=0}^{T-1}\E\mu_M(\vx_{t+1})&\leq \frac{\sqrt{M}F_0}{T} - \frac{M^{3/2}}{72} \E \frac{1}{T}\sum_{t=0}^{T-1}r_{t}^3 + 4 \frac{1}{T}\sum_{t=0}^{T-1}\E\|\veps_t\|^{3/2} + \frac{73 }{M^{3/2}}\frac{1}{T}\sum_{t=0}^{T-1}\E\| \mSigma_t \|^3\\
    &\leq \frac{\sqrt{M}F_0}{T} +  \frac{2}{\alpha T} \sigma_{g,0}^{3/2} +2   \alpha^{3/4}\sigma_g^{3/2} + 2\frac{1}{\gamma^{3/2}\alpha^{3/2}} \sigma_h^3  + \frac{73 }{M^{3/2}}\sigma_h^3\\
    &=\cO\Big( \frac{\sqrt{M}F_0}{T} +  \frac{M}{LT} \sigma_{g,0}^{3/2}  + \frac{ \sigma_g^{9/10}\sigma_h^{3/5}}{M^{3/10}}  + \frac{L^{3/4}}{M^{3/4}} \sigma_g^{3/2} +  \frac{1 }{M^{3/2}}\sigma_h^3\Big).
\end{align*}
We see that we have a term $ \frac{ \sigma_g^{9/10}\sigma_h^{3/5}}{M^{3/10}}$ that is much slower (for large $M$) than the $\frac{L^{3/8}}{M^{3/8}} \sigma_g^{3/2}$ we had before.

By sampling an initial batch $T/100$ larger than other batches, we ensure $\sigma_{g,0}^{3/2} = \cO(\frac{\sigma_{g}^{3/2}}{T^{3/4}}).$

We can optimize the right-hand side over $M$ as before, but we can already see that the obtained rate will worsen. Note also that the Hessian term $\frac{1 }{M^{3/2}}\sigma_h^3$ is smaller than the one we had before (since we needed to set $\beta=1$).

To be specific, by choosing $M=\max(3157 L, \frac{\sigma_g^{9/8}\sigma_h^{3/4} T^{5/4}}{F_0^{5/4}}, \frac{\sigma_g^{6/5}L^{3/5} T^{4/5}}{F_0^{4/5}},\frac{\sigma_h^2 T^{1/2}}{F_0^{1/2}})$, we get:
\begin{align*}
    \frac{1}{T}\sum_{t=0}^{T-1}\E\mu_M(\vx_{t+1})&=\cO\Big( \frac{\sqrt{L}F_0}{T} + \frac{\sigma_g^{9/16}\sigma_h^{3/8} F_0^{3/8}}{T^{3/8}} + \frac{L^{3/10}\sigma_g^{3/5} F_0^{3/5}}{T^{3/5}} +  \frac{\sigma_h^{3/4} F_0^{3/4} }{T^{3/4}}+  \frac{M}{LT^{7/4}} \sigma_{g}^{3/2} \Big)
\end{align*}
The term $\frac{M}{LT^{7/4}} \sigma_{g}^{3/2}$ leads to smaller terms compared to the ones we already have.

\section{Analyzing the HB momentum for gradient and Hessian.}\label{sec:HB}
We consider the Heavy-Ball (HB) momentum defined by the recursive relation:
\begin{equation}\label{HB}
    \vg_t = (1 - \alpha) \vg_{t-1} + \alpha \nabla_{\xi_t} f(\vx_t)\;.
\end{equation}
To analyze this momentum, we need the following (additional) assumption:
\begin{framedassumption}\label{smoothness}
    The gradients of the function $f$ are $L_g-$Lipschitz for some constant $L_g\geq 0$, meaning that for all $\vx,\vy$ we have $$\|\nabla f(\vx) - \nabla f(\vy)\|\leq L_g\|\vx - \vy\|\;.$$ 
\end{framedassumption}

We prove the following Lemma:
\begin{framedlemma}
    Under Assumption~\ref{smoothness}, for any $\gamma>0$:
    $$\frac{1}{T}\sum_{t=0}^{T-1} \E\|\veps_t\|^{3/2} \leq \sqrt{3}\Bigg(\frac{1}{\alpha T}\sigma_{g,0}^{3/2} + \alpha^{3/4}\sigma_{g}^{3/2} + (1-\alpha)^{3/2} \frac{L_g^{3/2}}{\alpha^{3/2}\gamma^{3/2}} + (1-\alpha)^{3/2} \frac{L_g^{3/2}\gamma^{3/2}}{\alpha^{3/2}}\frac{1}{T}\sum_{i=0}^{t-1} \E r_{t}^{3} \Bigg)$$
\end{framedlemma}
\begin{proof}

For this momentum, we have:

\begin{align*}
   \veps_t &= (1-\alpha) \vg_{t-1} + \alpha \nabla f_{\xi_t}(\vx_t) - \nabla f(\vx_t) \\
    &= (1-\alpha) \Big(\vg_{t-1} - \nabla f(\vx_t)\Big) + \alpha \underbrace{\Big(\nabla f_{\xi_t}(\vx_t) - \nabla f(\vx_t)\Big)}_{:=\vn_t} \\
     &= (1-\alpha) \Big(\veps_{t-1}+ \underbrace{\nabla f(\vx_{t-1}) - \nabla f(\vx_{t})}_{:=\vz(\vx_{t-1},\vx_{t})}\Big) + \alpha \vn_t  \\
     &= (1-\alpha)\Big(\veps_{t-1}+ \vz(\vx_{t-1},\vx_{t})\Big) + \alpha \vn_t, 
     \intertext{where we defined $\vz(\va,\vb):=\nabla f(\va) - \nabla f(\vb)$. According to Assumption~\ref{smoothness}: $\|\vz(\va,\vb)\|\leq L_g \|\va - \vb\|$.}
\end{align*}

Thus
$$\veps_t = (1 - \alpha)^t\veps_0 + (1-\alpha)\sum_{i=0}^{t-1} (1-\alpha)^i\vz(\vx_{t-i-1},\vx_{t-i}) + \alpha \sum_{i=0}^{t-1} (1-\alpha)^i\vn_{t-i}\;,\; t\geq 1.$$

And 
\begin{align*}
    \|\veps_t\|&\leq (1 - \alpha)^t\|\veps_0\| + (1-\alpha)L_g \sum_{i=0}^{t-1} (1-\alpha)^i r_{t-i}  + \alpha\|\sum_{i=0}^{t-1} (1-\alpha)^i\vn_{t-i}\|.
    \intertext{Now we raise this to the power $3/2$ and using the convexity of $x\mapsto x^{3/2}$, we get:}
    \|\veps_t\|^{3/2}&\leq \sqrt{3}\Bigg((1 - \alpha)^{3t/2}\|\veps_0\|^{3/2} + (1-\alpha)^{3/2} L_g^{3/2}\Big(\sum_{i=0}^{t-1} (1-\alpha)^i r_{t-i}\Big)^{3/2} + \alpha^{3/2}\|\sum_{i=0}^{t-1} (1-\alpha)^i\vn_{t-i}\|^{3/2}\Bigg)\\
    &\leq \sqrt{3}\Bigg((1 - \alpha)^{3t/2}\|\veps_0\|^{3/2} + (1-\alpha)^{3/2} L_g^{3/2}\Big(\sum_{i=0}^{t-1} (1-\alpha)^i r_{t-i}\Big)^{3/2} + \alpha^{3/2}\|\sum_{i=0}^{t-1} (1-\alpha)^i\vn_{t-i}\|^{3/2}\Bigg)\\
    &\leq \sqrt{3}\Bigg((1 - \alpha)^{3t/2}\|\veps_0\|^{3/2} + (1-\alpha)^{3/2} \frac{L_g^{3/2}}{\sqrt{\alpha}}\sum_{i=0}^{t-1} (1-\alpha)^i r_{t-i}^{3/2} + \alpha^{3/2}\|\sum_{i=0}^{t-1} (1-\alpha)^i\vn_{t-i}\|^{3/2}\Bigg).
    \intertext{Now we take the expectation, and in the same lines as for the (IT) momentum, we get:}
    \E\|\veps_t\|^{3/2} &\leq \sqrt{3}\Bigg((1 - \alpha)^{3t/2}\sigma_{g,0}^{3/2} + (1-\alpha)^{3/2} \frac{L_g^{3/2}}{\sqrt{\alpha}}\sum_{i=0}^{t-1} (1-\alpha)^i \E r_{t-i}^{3/2} + \alpha^{3/4}\sigma_{g}^{3/2}\Bigg).
    \intertext{We take now the average over iterations:}
    \frac{1}{T}\sum_{t=0}^{T-1}\E\|\veps_t\|^{3/2} &\leq \sqrt{3}\Bigg(\frac{1}{\alpha T}\sigma_{g,0}^{3/2} + (1-\alpha)^{3/2} \frac{L_g^{3/2}}{\alpha^{3/2}}\frac{1}{T}\sum_{t=0}^{t-1}  \E r_{t}^{3/2} + \alpha^{3/4}\sigma_{g}^{3/2}\Bigg).
    \intertext{Then we use Young's inequality to transform the $r^{3/2}$ terms to $r^3$ terms, we get for any $\gamma>0$:}
    \frac{1}{T}\sum_{t=0}^{T-1} \E\|\veps_t\|^{3/2} &\leq \sqrt{3}\Bigg(\frac{1}{\alpha T}\sigma_{g,0}^{3/2} + \alpha^{3/4}\sigma_{g}^{3/2} + (1-\alpha)^{3/2} \frac{L_g^{3/2}}{\alpha^{3/2}\gamma^{3/2}} + (1-\alpha)^{3/2} \frac{L_g^{3/2}\gamma^{3/2}}{\alpha^{3/2}}\frac{1}{T}\sum_{i=0}^{t-1} \E r_{t}^{3} \Bigg).
\end{align*}
This finishes the proof.
\end{proof}

\textbf{Convergence analysis.}

\begin{framedtheorem}
Using Cubic Newton with the gradient and Hessian updates in \eqref{HB} and \eqref{Hmom} respectively with $\alpha=\min(1,\frac{18 L_g^{4/5}}{M^{2/5}\sigma_g^{2/5}}), \beta=1$ and under the additional Assumption~\ref{smoothness}, is such that  for any $M\geq L$ we have
    $$\frac{1}{1008} \E\mu_M(\vx_{out}^T)\leq \frac{\sqrt{M}\Tilde{F_0}_M}{T} +  71\frac{L_g^{3/5}\sigma_g^{6/5}}{M^{3/10}} +  \frac{\Tilde{\sigma}_h^3}{M^{3/2}}$$
    where $\Tilde{F_0}_M = F_0 + \frac{12\sigma_g^{2/5}}{18 L_g^{4/5} M^{1/10}} \sigma_{g,0}^{3/2}\leq \Tilde{F_0}_L:=F_0 + \frac{12\sigma_g^{2/5}}{18 L_g^{4/5} L^{1/10}} \sigma_{g,0}^{3/2}$ and $\Tilde{\sigma}_h^3 = 12^{15/4} L_g^3 + 73 \sigma_h^3$., the last inequality is obtained because $M\leq L$.
\end{framedtheorem}
\begin{proof}

We consider $\beta = 1$ (other values do not improve the rate),
Applying Theorem~\ref{ThOneStep} \begin{align*}
    \frac{1}{1008} \frac{1}{T}\sum_{t=0}^{T-1}\E\mu_M(\vx_{t+1})&\leq \frac{\sqrt{M}F_0}{T} - \frac{M^{3/2}}{72} \E \frac{1}{T}\sum_{t=0}^{T-1}r_{t}^3 + 4 \frac{1}{T}\sum_{t=0}^{T-1}\E\|\veps_t\|^{3/2} + \frac{73 }{M^{3/2}}\frac{1}{T}\sum_{t=0}^{T-1}\E\| \mSigma_t \|^3\\
    &\leq \frac{\sqrt{M}F_0}{T} +  \frac{12}{\alpha T} \sigma_{g,0}^{3/2} + 8\alpha^{3/4}\sigma_{g}^{3/2} + 8 (1-\alpha)^{3/2} \frac{L_g^{3/2}}{\alpha^{3/2}\gamma^{3/2}}  + \frac{73 }{M^{3/2}} \sigma_h^3\\& + \Big((1-\alpha)^{3/2} \frac{L_g^{3/2}\gamma^{3/2}}{\alpha^{3/2}} - \frac{M^{3/2}}{72}\Big)\E\frac{1}{T}\sum_{t=0}^{T-1}r_{t}^3.
    \intertext{We choose $\alpha$ such that $ (1-\alpha)^{3/2} \frac{L_g^{3/2}\gamma^{3/2}}{\alpha^{3/2}} - \frac{M^{3/2}}{72}\leq 0$, for example $\alpha = 18 \frac{L_g\gamma}{M}$, thus  }
    \frac{1}{1008} \frac{1}{T}\sum_{t=0}^{T-1}\E\mu_M(\vx_{t+1})&\leq \frac{\sqrt{M}F_0}{T} +  \frac{12}{\alpha T} \sigma_{g,0}^{3/2} +  70 (\frac{L_g\gamma}{M})^{3/4}\sigma_{g}^{3/2} + 0.1 \frac{M^{3/2}}{\gamma^{3}}  + \frac{73 }{M^{3/2}} \sigma_h^3.
    \intertext{We choose $\gamma \in \argmin_{\gamma \in (0,\frac{M}{18 L_g}] }70 (\frac{L_g\gamma}{M})^{3/4}\sigma_{g}^{3/2} + 0.1 \frac{M^{3/2}}{\gamma^{3}} $, we find $\gamma = \min(\frac{M}{18 L_g}, \frac{M^{3/5}}{L_g^{1/5}\sigma_g^{2/5}})$. Therefore, $\alpha=\min(1,\frac{18 L_g^{4/5}}{M^{2/5}\sigma_g^{2/5}})$ and we get:}
    \frac{1}{1008} \frac{1}{T}\sum_{t=0}^{T-1}\E\mu_M(\vx_{t+1})&\leq \frac{\sqrt{M}F_0}{T} +  \frac{12M^{2/5}\sigma_g^{2/5}}{18 L_g^{4/5} T} \sigma_{g,0}^{3/2} +  71\frac{L_g^{3/5}\sigma_g^{6/5}}{M^{3/10}} + 12^{15/4} \frac{L_g^3}{M^{3/2}}  + \frac{73 }{M^{3/2}} \sigma_h^3.
    \intertext{The term $\frac{L_g^3}{M^{3/2}}$ comes from the constraint $\gamma\leq \frac{M}{18 L_g}$ being activated which means $\sigma_g^{3/2}\leq 12^{15/4} \frac{L_g^3}{M^{3/2}} $. This term limits any possible improvement of the Hessian term (since it has the same power of $M$), hence why we chose $\beta=1$ from the start.}
    \intertext{Finally, we obtain}
    \frac{1}{1008} \frac{1}{T}\sum_{t=0}^{T-1}\E\mu_M(\vx_{t+1})&\leq \frac{\sqrt{M}\Tilde{F_0}_M}{T} +  71\frac{L_g^{3/5}\sigma_g^{6/5}}{M^{3/10}} +  \frac{12^{15/4} L_g^3 + 73 \sigma_h^3 }{M^{3/2}},
    \intertext{where $\Tilde{F_0}_M = F_0 + \frac{12\sigma_g^{2/5}}{18 L_g^{4/5} M^{1/10}} \sigma_{g,0}^{3/2}\leq \Tilde{F_0}_L:=F_0 + \frac{12\sigma_g^{2/5}}{18 L_g^{4/5} L^{1/10}} \sigma_{g,0}^{3/2}$, the last inequality is obtained because $M\leq L$.}
\end{align*}

\end{proof}
 Choosing $M = \max(L,\frac{L_g^{3/4} \sigma_g^{3/2} T^{5/4}}{\Tilde{F_0}_L^{5/4}}, \frac{\Tilde{\sigma}_h^2 T^{1/2}}{\Tilde{F_0}_L^{1/2}}$), we get:
 $$\frac{1}{1008}\frac{1}{T}\sum_{t=0}^{T-1}\E\mu_M(\vx_{t+1}) \leq \frac{\sqrt{L}\Tilde{F_0}_L}{T} + 72 \Big(\frac{L_g\Tilde{F_0}_L\sigma_g^2}{T}\Big)^{3/8} + \Big(\frac{\Tilde{F_0}_L\Tilde{\sigma}_h}{T}\Big)^{3/4}$$

%###################################################
\section{Analyzing the MVR momentum for the gradient.}\label{sec:MVR}

We consider the Heavy-Ball (MVR) momentum defined by the recursive relation:
\begin{equation}\label{MVR}
    \vg_t = (1 - \alpha) \Big(\vg_{t-1} + \nabla_{\xi_t} f(\vx_t) - \nabla_{\xi_t} f(\vx_{t-1})\Big) + \alpha \nabla_{\xi_t} f(\vx_t)\;.
\end{equation}
To analyze this momentum, we need a stronger first-order smoothness assumption:
\begin{framedassumption}\label{Gensmoothness}
    There exists $L_g\geq 0$ such that for all $\xi$, the gradients of the function $f_\xi$ are $L_g-$Lipschitz, meaning that for all $\xi,\vx,\vy$ we have $$\|\nabla f_\xi(\vx) - \nabla f_\xi(\vy)\|\leq L_g\|\vx - \vy\|\;.$$ 
\end{framedassumption}

We prove the following Lemma:
\begin{framedlemma}
For the gradient update defined in \eqref{MVR}, we have
    for any $\gamma>0$:
    $$\frac{1}{T}\sum_{t=0}^{T-1} \E\|\veps_t\|^{3/2} \leq \frac{1}{\alpha T}\sigma_{g,0}^{3/2} + \alpha^{3/2} \sigma_g^{3/2} + \frac{8 (1- \alpha)^{3/2} L_g^{3/2}}{\alpha\gamma^{3/2}} + \frac{8 (1- \alpha)^{3/2} L_g^{3/2}\gamma^{3/2}}{\alpha} \frac{1}{T}\sum_{i=0}^{t-1} r_{t}^{3} $$
\end{framedlemma}
\begin{proof}

For this momentum, we have:

\begin{align*}
   \veps_t &= (1-\alpha) \Big(\vg_{t-1}+ \nabla_{\xi_t} f(\vx_t) - \nabla_{\xi_t} f(\vx_{t-1})\Big) + \alpha \nabla f_{\xi_t}(\vx_t) - \nabla f(\vx_t) \\
    &= (1-\alpha) \Big(\vg_{t-1} + \nabla_{\xi_t} f(\vx_t) - \nabla_{\xi_t} f(\vx_{t-1}) - \nabla f(\vx_t)\Big) + \alpha \underbrace{\Big(\nabla f_{\xi_t}(\vx_t) - \nabla f(\vx_t)\Big)}_{:=\vn_t} \\
     &= (1-\alpha) \Big(\veps_{t-1}+ \underbrace{\nabla f(\vx_{t-1}) - \nabla f(\vx_{t})+ \nabla_{\xi_t} f(\vx_t) - \nabla_{\xi_t} f(\vx_{t-1})}_{:=\vz_{\xi_t}(\vx_{t-1},\vx_{t})}\Big) + \alpha \vn_t  \\
     &= (1-\alpha)\Big(\veps_{t-1}+ \vz_{\xi_t}(\vx_{t-1},\vx_{t})\Big) + \alpha \vn_t, 
     \intertext{where we defined $\vz_{\xi}(\va,\vb):=\nabla f(\va) - \nabla f(\vb) + \nabla f_{\xi}(\vb) - \nabla f_{\xi}(\va)$. According to Assumption~\ref{Gensmoothness}: $\|\vz_{\xi}(\va,\vb)\|\leq 2L_g \|\va - \vb\|$.}
\end{align*}
Note that knowing $\vx_{t-1},\vx_t$, the random variable $(1- \alpha)\vz_{\xi_t}(\vx_{t-1},\vx_{t}) + \alpha \vn_t$ is centred and independent from $\veps_{t-1}$.

Thus
\begin{align*}
     \E\|\veps_t\|^2 &= (1 - \alpha)^2 \E\|\veps_{t-1}\|^2 + \E\|(1- \alpha)\vz_{\xi_t}(\vx_{t-1},\vx_{t}) + \alpha \vn_t\|^2\\
     &\leq (1 - \alpha)^2 \E\|\veps_{t-1}\|^2 + 8 (1- \alpha)^2 L_g^2 r_{t}^2 + 2 \alpha^2 \sigma_g^2\\
     &\leq (1 - \alpha)^{2t}\|\veps_0\|^2 + 8 (1- \alpha)^2 L_g^2 \sum_{i=0}^{t-1} (1-\alpha)^{2i} r_{t-i}^2  + \alpha^2\sum_{i=0}^{t-1} (1-\alpha)^{2i} \sigma_g^2\\
     &\leq (1 - \alpha)^{2t}\sigma_{g,0}^2 + 8 (1- \alpha)^2 L_g^2 \sum_{i=0}^{t-1} (1-\alpha)^{2i} r_{t-i}^2  + \alpha \sigma_g^2.
     \intertext{Using Jensen's inequality, we get:}
     \E\|\veps_t\|^{3/2}&\leq (1 - \alpha)^{3t/2}\sigma_{g,0}^{3/2} + 8 (1- \alpha)^{3/2} L_g^{3/2} \Big(\sum_{i=0}^{t-1} (1-\alpha)^{2i} r_{t-i}^2\Big)^{3/4}  + \alpha^{3/2} \sigma_g^{3/2}\\
     &\leq (1 - \alpha)^{3t/2}\sigma_{g,0}^{3/2} + 8 (1- \alpha)^{3/2} L_g^{3/2} \sum_{i=0}^{t-1} (1-\alpha)^{3i/2} r_{t-i}^{3/2}  + \alpha^{3/2} \sigma_g^{3/2}.\\
     \intertext{We take the average and get:}
     \frac{1}{T}\sum_{t=0}^{T-1}\E\|\veps_t\|^{3/2}&\leq \frac{1}{\alpha T}\sigma_{g,0}^{3/2} + \alpha^{3/2} \sigma_g^{3/2} + \frac{8 (1- \alpha)^{3/2} L_g^{3/2}}{\alpha} \frac{1}{T}\sum_{i=0}^{t-1} r_{t}^{3/2}\\  
     &\leq \frac{1}{\alpha T}\sigma_{g,0}^{3/2} + \alpha^{3/2} \sigma_g^{3/2} + \frac{8 (1- \alpha)^{3/2} L_g^{3/2}}{\alpha\gamma^{3/2}} + \frac{8 (1- \alpha)^{3/2} L_g^{3/2}\gamma^{3/2}}{\alpha} \frac{1}{T}\sum_{i=0}^{t-1} r_{t}^{3},
     \intertext{for any $\gamma>0$. This finishes the proof.}
\end{align*}
\end{proof}
\begin{framedtheorem}
Using Cubic Newton with the gradient and Hessian updates in \eqref{MVR} and \eqref{Hmom} respectively with $\alpha=\min(1,\frac{18 L_g^{4/5}}{M^{2/5}\sigma_g^{2/5}}), \beta=1$ and under the additional Assumption~\ref{Gensmoothness}, is such that  for any $M\geq L$ we have
    $$\frac{1}{1008} \E\mu_M(\vx_{out}^T)\leq\frac{\sqrt{M}\Tilde{F_0}_M}{T} +  \frac{M^{6/11}\sigma_g^{6/11}}{192 L_g^{12/11} T} \sigma_{g,0}^{3/2} +327  \frac{L_g^{9/11}\sigma_g^{12/11}}{M^{9/22}} +  \frac{\Tilde{\sigma}_h^3 }{M^{3/2}} ,$$
    where $\Tilde{F_0}_M = F_0 + \frac{M^{1/22}\sigma_g^{6/11}}{192 L_g^{12/11}} \sigma_{g,0}^{3/2} $. and $\Tilde{\sigma}_h^3 = (72\times 32)^{11/6} L_g^3 +73 \sigma_h^3$.
\end{framedtheorem}
\begin{proof}

We consider $\beta = 1$ (other values do not improve the rate),
Applying Theorem~\ref{ThOneStep} \begin{align*}
    \frac{1}{1008} \frac{1}{T}\sum_{t=0}^{T-1}\E\mu_M(\vx_{t+1})&\leq \frac{\sqrt{M}F_0}{T} - \frac{M^{3/2}}{72} \E \frac{1}{T}\sum_{t=0}^{T-1}r_{t}^3 + 4 \frac{1}{T}\sum_{t=0}^{T-1}\E\|\veps_t\|^{3/2} + \frac{73 }{M^{3/2}}\frac{1}{T}\sum_{t=0}^{T-1}\E\| \mSigma_t \|^3\\
    &\leq \frac{\sqrt{M}F_0}{T} +  \frac{12}{\alpha T} \sigma_{g,0}^{3/2} + 32 \alpha^{3/4}\sigma_{g}^{3/2} + 32 (1-\alpha)^{3/2} \frac{L_g^{3/2}}{\alpha\gamma^{3/2}}  + \frac{73 }{M^{3/2}} \sigma_h^3\\& + \Big(32 (1-\alpha)^{3/2} \frac{L_g^{3/2}\gamma^{3/2}}{\alpha} - \frac{M^{3/2}}{72}\Big)\E\frac{1}{T}\sum_{t=0}^{T-1}r_{t}^3.
    \intertext{We choose $\alpha$ such that $ 32  \frac{L_g^{3/2}\gamma^{3/2}}{\alpha} - \frac{M^{3/2}}{72}\leq 0$, for example $\alpha = 72\times 32 \big(\frac{L_g\gamma}{M}\big)^{3/2}$, thus  }
    \frac{1}{1008} \frac{1}{T}\sum_{t=0}^{T-1}\E\mu_M(\vx_{t+1})&\leq \frac{\sqrt{M}F_0}{T} +  \frac{12}{\alpha T} \sigma_{g,0}^{3/2} +  326 (\frac{L_g\gamma}{M})^{9/8}\sigma_{g}^{3/2} + \frac{M^{3/2}}{72\gamma^{3}}  + \frac{73 }{M^{3/2}} \sigma_h^3.
    \intertext{We choose $\gamma \in \argmin_{\gamma \in (0,\frac{M}{172 L_g}] } 326 (\frac{L_g\gamma}{M})^{9/8}\sigma_{g}^{3/2} + \frac{M^{3/2}}{72\gamma^{3}} $, we find $\gamma = \min(\frac{M}{172 L_g}, \frac{M^{7/11}}{L_g^{3/11}\sigma_g^{4/11}})$, meaning $\alpha=\min\Big(1,72\times 32 (\frac{L^2_g}{M\sigma_g})^{6/11}\Big)$ and get:}
    \frac{1}{1008} \frac{1}{T}\sum_{t=0}^{T-1}\E\mu_M(\vx_{t+1})&\leq \frac{\sqrt{M}F_0}{T} +  \frac{M^{6/11}\sigma_g^{6/11}}{192 L_g^{12/11} T} \sigma_{g,0}^{3/2} +327  \frac{L_g^{9/11}\sigma_g^{12/11}}{M^{9/22}} + (72\times 32)^{11/6} \frac{L_g^3}{M^{3/2}}  + \frac{73 }{M^{3/2}} \sigma_h^3.
    \intertext{The term $\frac{L_g^3}{M^{3/2}}$ comes from the constraint $\gamma\leq \frac{M}{172 L_g}$ being activated which means $\sigma_g^{3/2}\leq (72\times 32)^{11/6} \frac{L_g^3}{M^{3/2}} $. This term limits any possible improvement of the Hessian term (since it has the same power of $M$), hence why we chose $\beta=1$ from the start.}
    \intertext{Finally, we get}
    \frac{1}{1008} \frac{1}{T}\sum_{t=0}^{T-1}\E\mu_M(\vx_{t+1})&\leq\frac{\sqrt{M}\Tilde{F_0}_M}{T} +  \frac{M^{6/11}\sigma_g^{6/11}}{192 L_g^{12/11} T} \sigma_{g,0}^{3/2} +327  \frac{L_g^{9/11}\sigma_g^{12/11}}{M^{9/22}} +  \frac{\Tilde{\sigma}_h^3 }{M^{3/2}} ,
    \intertext{where $\Tilde{F_0}_M = F_0 + \frac{M^{1/22}\sigma_g^{6/11}}{192 L_g^{12/11} } \sigma_{g,0}^{3/2} $, and $\Tilde{\sigma}_h^3 = (72\times 32)^{11/6} L_g^3 +73 \sigma_h^3$.}
\end{align*}
\end{proof}
Note here that the power of $M$ of the gradient term (meaning $\frac{L_g^{9/11}\sigma_g^{12/11}}{M^{9/22}}$) is better than the one we had using the IT momentum in \eqref{Gradmom} ($\frac{9}{22} - \frac{3}{8} = \frac{6}{176} > 0$), thus we should expect to have a better dependence on $T$ for this term.

By sampling a large first batch (for example $b_g = T/100$ for the gradient such that we can ensure $\sigma_{g,0}\leq \frac{10\sigma_g}{\sqrt{T}}$, and

choosing $M = \max(L,\frac{L_g^{9/10} \sigma_g^{6/5} T^{11/10}}{F_0^{11/10}}, \frac{\Tilde{\sigma}_h^2 T^{1/2}}{F_0^{1/2}}$), we get:
\begin{align*}
    \frac{1}{1008}\frac{1}{T}\sum_{t=0}^{T-1}\E\mu_M(\vx_{t+1}) &\leq \frac{\sqrt{M}F_0}{T} +327  \frac{L_g^{9/11}\sigma_g^{12/11}}{M^{9/22}} +  \frac{\Tilde{\sigma}_h^3}{M^{3/2}} +  \frac{M^{6/11}\sigma_g^{45/22}}{192 L_g^{12/11} T^{7/4}}\\&\leq \frac{\sqrt{L}F_0}{T} + 328 \Big(\frac{L_g F_0}{T}\Big)^{9/20}\sigma_g^{6/10} + \Big(\frac{F_0\Tilde{\sigma}_h}{T}\Big)^{3/4}.
\end{align*}

Note that $\frac{1}{T^{9/20}} \leq \frac{1}{T^{3/7}}$, which means that the MVR momentum slightly improves the gradient rate of stochastic cubic Newton compared to the IT momentum, but this comes at the cost of a slower rate for the Hessian term.

\end{document}

%% file: math_commands.tex
\usepackage{amsfonts,bm}

\usepackage{enumitem}

\newcommand{\argmin}{{\rm arg\,min}}

\def\beq{\begin{equation}}
\def\eeq{\end{equation}}
\def\ba{\begin{array}}
\def\ea{\end{array}}
\def\beann{\begin{eqnarray*}}
\def\eeann{\end{eqnarray*}}
\def\bea{\begin{eqnarray}}
\def\eea{\end{eqnarray}}

\def\BT{\begin{theorem}}
\def\ET{\end{theorem}}
\def\BL{\begin{lemma}}
\def\EL{\end{lemma}}
\def\BC{\begin{corollary}}
\def\EC{\end{corollary}}
\def\BE{\begin{example}}
\def\EE{\end{example}}
\def\BD{\begin{definition}}
\def\ED{\end{definition}}
\def\BR{\begin{remark}}
\def\ER{\end{remark}}
\def\BAS{\begin{assumption}}
\def\EAS{\end{assumption}}
\def\BI{\begin{itemize}}
\def\EI{\end{itemize}}
\def\BP{\begin{proposition}}
\def\EP{\end{proposition}}

\def\BMP{\begin{minipage}{9.5cm}}
\def\EMP{\end{minipage}}
\def\MPT{\begin{minipage}{11.5cm}}
\def\EPT{\end{minipage}}

\def\la{\langle}
\def\ra{\rangle}

\def\1{\bm{1}}

\def\vzero{{\bm{0}}}

\def\veps{{\bm{\varepsilon}}}
\def\vtheta{{\bm{\theta}}}
\def\va{{\bm{a}}}
\def\vb{{\bm{b}}}

\def\vg{{\bm{g}}}

\def\vn{{\bm{n}}}

\def\vx{{\bm{x}}}
\def\vy{{\bm{y}}}
\def\vz{{\bm{z}}}

\def\mA{{\bm{A}}}

\def\mH{{\bm{H}}}
\def\mI{{\bm{I}}}

\def\mN{{\bm{N}}}

\def\mZ{{\bm{Z}}}

\def\mSigma{{\bm{\Sigma}}}

\DeclareMathAlphabet{\mathsfit}{\encodingdefault}{\sfdefault}{m}{sl}
\SetMathAlphabet{\mathsfit}{bold}{\encodingdefault}{\sfdefault}{bx}{n}

\def\cO{\mathcal{O}}
\def\cB{\mathcal{B}}

\newcommand{\E}{\mathbb{E}}

\newcommand{\R}{\mathbb{R}}

\renewcommand{\epsilon}{\varepsilon}